\documentclass[12pt]{iopart}
\usepackage{iopams}
\usepackage{graphicx}
\usepackage{showlabels}
\usepackage[usenames,dvipsnames]{color}
\usepackage[normalem]{ulem}
\usepackage{multirow}
\input{epsf}
\newcommand{\hyp}{\,_2F_1}
\newcommand{\FT}{\widehat{P}_z}

\begin{document}

\title[Non-homogeneous random walks with resetting]{Non-homogeneous random walks with stochastic resetting: an application to the Gillis model}

\author{Mattia Radice}
\address{Max Planck Institute for the Physics of Complex Systems, N\"{o}thnitzer Stra\ss e 38, 01187 Dresden, Germany}
\eads{\mailto{mradice@pks.mpg.de}}
\vspace{10pt}

\begin{abstract}
We consider the problem of the first passage time to the origin of a spatially non-homogeneous random walk with a position-dependent drift, known as the \textit{Gillis random walk}, in the presence of resetting. The walk starts from an initial site $ x_0 $ and, with fixed probability $ r $, at each step may be relocated to a given site $ x_r $. From a general perspective, we first derive a series of results regarding the first and the second moment of the first hitting time distribution, valid for a wide class of processes, including random walks lacking the property of translational invariance; we then apply these results to the specific model. When resetting is not applied, by tuning the value of a parameter which defines the transition probability of the process, denoted by $ \epsilon $, the recurrence properties of the walk are changed, and we can observe: a transient walk, a null-recurrent walk, or a positive-recurrent walk. When the resetting mechanism is switched on, we study quantitatively in all regimes the improvement of the search efficiency. In particular, in every case resetting allows the system to reach the target with probability one and, on average, in a finite time. If the reset-free system is in the transient or null-recurrent regime, this makes resetting always advantageous and moreover, it assures the existence of an optimal resetting probability $ r^* $ which minimizes the mean first hitting time. Instead, when the system is positive-recurrent, the introduction of resetting is not necessarily beneficial. We explain that in this case there exists a threshold $ r_{\mathrm{th}} $ for the resetting probability $ r $, above which the resetting mechanism yields a larger mean first hitting time with respect to the reset-free system. We provide a study of $ r_{\mathrm{th}} $, which can be zero for some values of the system parameters, meaning that resetting can not be beneficial in those situations. All the theoretical findings are corroborated with numerical simulations.
\end{abstract}

\vspace{2pc}

\section{Introduction}
Stochastic resetting is a simple mechanism by which from time to time a random process is restored to a given condition, usually corresponding to the initial state \cite{EvaMaj-2011}. Following each resetting event the dynamics is started anew, with loss of memory of the previous history. Despite the simplicity of the mechanism, the repeated returns to a given location produce non-trivial effects on the statistical properties of the dynamics. For instance, if we consider the paradigmatic example of Brownian motion, it is known that resetting favours the emergence of a non-equilibrium steady state \cite{EvaMaj-2014}, which arises because the diffusive spreading of the probability density function is suppressed and the particle is only allowed to diffuse between to consecutive restarts. A second striking fact regards the time needed to reach a target for the first time, which is often investigated to test the efficiency of a search process. While for a Brownian motion without restarts the expected time to reach a target is infinite, the introduction of resetting eventually cuts off the contribution of long diffusive excursions away from the target, thus enhancing the search efficiency to the point where the mean first passage time attains a finite value. Not surprisingly, this property has been thoroughly investigated in the literature \cite{EvaMaj-2011,Reu-2016,PalKunEva-2016,PalReu-2017,CheSok-2018,TalPalSek-2020,BesBovPet-2020,SinPal-2022,RAD-2022,EliReu} and it is still object of many studies, trying to extend it outside the realm of pure diffusion, e.g. for the telegraphic process and run-and-tumble motion \cite{EvaMaj-2018,Mas-2019,SanBasSab-2020,RAD-2021,TucGamMaj-2022}, continuous time random walks \cite{MonMasVil-2017,KusGod-2019} and Lévy flights \cite{KusMajSabSch-2014,KusGod-2015,StaWer-2021,ZhoXuDen-2021}. For comprehensive reviews, see \cite{EvaMajSch-2020,PalKosReu-2021}. It is worth mentioning that the usefulness of stochastic resetting has been widely recognized also in different contexts, such as quantum physics \cite{PerCarMag-2021,MagCarPer-2022,YinBar-2022}, chemistry \cite{ReuUrbKla-2014,RotReuUrb-2015}, biology \cite{RolLisSanGri-2016,RamMagSan-2020}, computer science \cite{LubSinZuc-1993,MonZec-2002} and economics \cite{StoSanKoc-2021,StoJolPal-2022, San-2022}.

Although the vast majority of the studies focused on situations where the time of the process is described by a continuous parameter, a reasonable amount of work has been dedicated to the analysis of random walks (RW) in discrete time. In \cite{MajSabSch-2015Res} the authors investigated a one-dimensional random walker on a discrete lattice with resetting to the position of the maximum. In \cite{MonVil-2016} it was considered the introduction of resetting in the so-called Sisyphus RW, a particular kind of unidirectional walk. Resetting with memory of the previous history was studied in for RW in discrete time and space \cite{BoyFalGiu-2019}. A general setup for the analysis of first passage RW processes was proposed in \cite{BonPal-2021}. Furthermore, stochastic resetting has been also considered for RWs on networks, with a particular attention on first passage times \cite{RiaBoyHer-2020,RiaBoyHer-2022}. More recently, RWs with stochastic resetting have been considered for the evaluation of extreme value statistics, regarding maxima, records and survival probability \cite{MajMouSabSch-2021,God-Luc-2022Rec,God-Luc-2022Sur}, and to revisit the classical problem of the number of distinct sites visited up to step $ n $ \cite{BirMorMaj-2022}.

In this paper we will present general results regarding the search efficiency of the process in presence of resetting, valid for a wide class of random walks on the integer lattice, which will be specified later. Our results involve the evaluation of the mean and the variance of the first hitting time to a target, with the aim of unveiling the relation with the recurrence properties of the reset-free system. To corroborate our general findings, we will consider a classical model of random walk on the one-dimensional lattice of integers known by the name of \textit{Gillis random walk} (GRW) \cite{Gill-bias}, defined by the following transition probability $ G_{xy} $ of jumping from site $ y $ to site $ x $:
\begin{equation}\label{eq:Trans_def}
	G_{xy}=\cases{
	\frac12\left(1-\frac{\epsilon}{y}\right)\delta_{x,y+1}+\frac12\left(1+\frac{\epsilon}{y}\right)\delta_{x,y-1}&if $y\neq0$\\
	\frac12\delta_{x,y+1}+\frac12\delta_{x,y-1}&if $y=0$,
	}
\end{equation}
where $ \delta_{xy} $ is the Kronecker delta, restricting the jumps only between nearest-neighbour sites, while $ \epsilon $ is a real parameter such that $ -1<\epsilon\leq1 $. As can be inspected from \eref{eq:Trans_def}, there is a tendency to move towards or away from the origin, depending on the sign of $ \epsilon $. Moreover, the transition probability depends explicitly on the current position. For these reasons, the GRW can be classified as a \textit{non-homogeneous centrally-biased random walk}. The peculiarity of this model is that is one the few kinds of random walks lacking translational invariance for which exact analytical results can be obtained \cite{Hug-I}. It recently attracted interest in the physical literature \cite{OPRA,PROA} due to the correspondence of the asymptotic properties of the walk to those of a Brownian particle in a logarithmic potential, even if the model already appeared in \cite{ChaHug-1988}, to map the problem of ion diffusion in a semi-infinite domain in the presence of a charged barrier. However, the GRW is mostly known in the mathematical literature as a paradigmatic example for the study of recurrence and first passage properties in $ d $ dimensions and in the presence of a bias \cite{Gill-bias,Nas-1959,Lam-1960,Lam-1963,HryMenWad-2013}.

The paper is organised as follows: first, in sections \ref{s:RecTra} and \ref{s:RW_with_res} we briefly review the notions of recurrence and transience, and present the formalism to treat analytically the first hitting problem in the presence of resetting; then, in sections \ref{s:Anal} and \ref{s:Anal2} we derive general results regarding the first and the second moment of the first hitting time distribution, also specifying the class of processes for which these are valid; in section \ref{s:Opt} we discuss under which conditions the problem may be optimised by choosing correctly the probability of reset; in section \ref{s:Results} we apply all findings to the Gillis model and show the agreement with our numerical simulations; finally in section \ref{s:Conc} we draw our conclusions.  

\section{Recurrence and transience}\label{s:RecTra}
In this section we recall the basic definitions of recurrence and transience for a random walk, and also specify the class of processes that we will consider in the paper. We say that a random walk is recurrent if with probability one it returns to the starting site an infinite number of times \cite{Gill-bias}. An useful criterion to test recurrence is thus based on the generating function of the probability of first return time at the origin. Let us denote with $ F(n) $ the probability that a walk starting at $ x_0=0 $ returns to the origin for the first time at step $ n $, and define the generating function
\begin{equation}\label{eq:GF_FRet}
	\widetilde{F}(z)=\sum_{n=1}^{\infty}F(n)z^n.
\end{equation}
The limit $ z\to 1^{-} $ of this generating function corresponds to the sum of all $ F(n) $, hence the random walk is recurrent if and only if
\begin{equation}
	\lim_{z\to1^{-}}\widetilde{F}(z)=1\qquad(\textrm{recurrent walk});
\end{equation}
otherwise we say that the random walk is transient, which means that there is a positive probability $ Q=1-R $ that the walk never returns to the starting point, where $ R $ is given by
\begin{equation}
	\lim_{z\to1^{-}}\widetilde{F}(z)=R<1\qquad(\textrm{transient walk}).
\end{equation}
It is important also to note that one can distinguish between \textit{null} and \textit{positive} recurrence: a null-recurrent walk returns with \textit{infinite} expected return time, while the return time for positive-recurrent walks has finite mean.

It turns out that when a random walk is recurrent, the general structure assumed by $ \widetilde{F}(z) $ is ruled by a classical result by Lamperti \cite{Lam}, which establishes a relation between $ \widetilde{F}(z) $ and the existence of a certain limiting distribution. To be more specific, the theorem is valid for a class of processes with the property that the state space can be split in two sets, say $\mathcal{S}_1$ and $ \mathcal{S}_2 $, separated by a special state $ s $, assumed as the initial condition, and furthermore: (i) the occupation of $ s $ is a recurrent event; (ii) if the process is in $ \mathcal{S}_1 $ at step $ n-1 $ and in $ \mathcal{S}_2 $ at step $ n+1 $ (or vice versa), then necessarily occupies $ s $ at step $ n $; under these condition, the existence for $ n\to\infty $ of the limiting distribution of the fraction of time spent in $ \mathcal{S}_1 $ implies that there is $ 0\leq\rho\leq1 $ and a slowly-varying function $ \mathcal{L}(x) $ such that
\begin{equation}\label{eq:F_form}
	\widetilde{F}(z)=1-(1-z)^{\rho}\mathcal{L}\left(\frac{1}{1-z}\right).
\end{equation}
We recall that a slowly-varying function is a continuous function, positive for large enough $ x $,  that for any $ t>0 $ satisfies
\begin{equation}\label{key}
	\lim_{x\to\infty}\frac{\mathcal{L}(tx)}{\mathcal{L}(x)}=1.
\end{equation}

The class of processes for which this result holds is quite general, and includes for example all Markov chains in which one state separates the rest in two sets \cite{Lam}. For example, considering the simple and symmetric random walk on $ \mathbb{Z} $ starting from $ x_0=0 $, then we can identify $ \mathcal{S}_1 $ with the set of positive integers and $ \mathcal{S}_2$ with the set of negative integers; the theorem states that if the distribution of the fraction of time spent in $ \mathcal{S}_1 $ exists, then $ \widetilde{F}(z) $ can be written as in \eref{eq:F_form}, and indeed we know \cite{Fell-I}
\begin{equation}\label{key}
	\widetilde{F}(z)=1-\sqrt{1-z^2}=1-\sqrt{1-z}\mathcal{L}\left(\frac{1}{1-z}\right),
\end{equation}
with
\begin{equation}\label{key}
	\mathcal{L}(x)=\sqrt{2-\frac1x},
\end{equation}
which shows that for this model $ \rho=\case12 $.

We remark that the values $ \rho=0 $ and $ \rho=1 $ are special, because they both define a “phase transition": walks with $ \rho=0 $ correspond to a transition between recurrent and transient processes, while those with $ \rho=1 $ separate null-recurrence and positive-recurrence \cite{Lam}. The consequences of this result on many statistical properties of the walk, as well as the evaluation of $ \rho $ and its general relation with other exponents describing the system, have been investigated in previous works \cite{OPRA,ROAP}. Now we are interested in the special form of $ \widetilde{F}(z) $ given by \eref{eq:F_form}, which will help us obtain straightforwardly useful information on the mean and the standard deviation of the first hitting time in presence of resetting.  

\section{Random walks with stochastic resetting}\label{s:RW_with_res}
In this section we are interested in the number of steps required to reach the point $ x $ for the first time for a random walker starting from $ x_0 $, when we also consider a stochastic resetting mechanism that, with probability $ r $, after each step can relocate the walk to a fixed position $ x_r $. This is sometimes called \textit{geometric restart} \cite{BonPal-2021}. In the following we will only consider walks that can perform jumps between nearest-neighbour sites. Let us first consider the simple case where $ x_0=x_r $, and let us denote by $ F_r(x,n|x_0) $ the probability of hitting $ x $ for the first time at step $ n $, in the presence of resetting which relocates the walker at the starting position $ x_0 $. This can be expressed as the sum of two contributions: the first comes from the walks that reach $ x $ in $ n $ steps with no resetting events. This can be written in terms of $ F(x,n|x_0) $, by which we denote the first hitting time probability of $ x $ in absence of resetting. The second is due to walks that are reset for the first time at step $ k<n $, not having hit $ x $ in the previous $ k-1 $ steps, and then reach the target in the remaining $ n-k $ steps. Hence by calling $ Q(x,n|x_0) $ the probability of never reaching $ x $ in $ n $ steps starting from $ x_0 $, namely the \textit{survival probability}, we can write
\begin{eqnarray}
	F_r(x,n|x_0)&=&(1-r)^nF(x,n|x_0)+\nonumber\\
	&&r\sum_{k=1}^{n}(1-r)^{k-1}Q(x,k-1|x_0)F_r(x,n-k|x_0).\label{eq:F_rEQ}
\end{eqnarray}
Note that an equivalent equation for the survival probability in presence of resetting, $ Q_r(x,n|x_0) $, already appeared, for instance, in \cite{KusMajSabSch-2014}.
We define the generating functions
\begin{eqnarray}
	\widetilde{F}_r(x,z|x_0)&=&\sum_{n=1}^{\infty}z^nF_r(x,n|x_0)\\
	\widetilde{F}(x,z|x_0)&=&\sum_{n=1}^{\infty}z^nF(x,n|x_0)\\
	\widetilde{Q}(x,z|x_0)&=&\sum_{n=0}^{\infty}z^nQ(x,n|x_0),
\end{eqnarray}
then multiply \eref{eq:F_rEQ} for $ z^n $ and sum over $ n=1,2,\dots $, obtaining
\begin{equation}\label{eq:F_rGenEQ}
	\widetilde{F}_r(x,z|x_0)=\widetilde{F}(x,\zeta|x_0)+zr\widetilde{Q}(x,\zeta|x_0)\widetilde{F}_r(x,z|x_0),
\end{equation}
where we have introduced the rescaled variable $ \zeta=(1-r)z $. It is more convenient to rewrite $ \widetilde{Q}(x,z|x_0) $ in terms of $ \widetilde{F}(x,z|x_0) $, which can be done by observing that the numbers of walks that reach $ x $ for the first time in $ n $ steps is the difference between those that have not reach it up to step $ n-1 $, and those that have not reach it up to step $ n $:
\begin{equation}\label{key}
	F(x,n|x_0)=Q(x,n-1|x_0)-Q(x,n|x_0).
\end{equation}
Hence one gets the following relation between the generating functions:
\begin{equation}\label{eq:QvsF}
	\widetilde{Q}(x,z|x_0)=\frac{1-\widetilde{F}(x,z|x_0)}{1-z}.
\end{equation}
By plugging this expression in \eref{eq:F_rGenEQ} and rearranging terms, we get
\begin{equation}\label{eq:FGF_r_x0}
	\widetilde{F}_r(x,z|x_0)=\frac{(1-\zeta)\widetilde{F}(x,\zeta|x_0)}{1-z+rz\widetilde{F}(x,\zeta|x_0)}.
\end{equation}

We can now easily generalize to the case of a generic resetting position $ x_r $. Obviously the difference with respect to the previous case comes only from those walks that experience at least one resetting, hence we have to modify only the second term at the right-hand side of \eref{eq:F_rEQ} to take into account that after the first reset, the process is reset not at $ x_0 $ but at $ x_r $, and that all future resets will relocate the particle at the same position $ x_r $. Therefore in this case we have
\begin{eqnarray}\label{key}
	F_r(x,n|x_0,x_r)&=&(1-r)^nF(x,n|x_0)+\nonumber\\
	&&r\sum_{k=1}^{n}(1-r)^{k-1}Q(x,k-1|x_0)F_r(x,n-k|x_r),
\end{eqnarray}
and by switching to generating functions:
\begin{equation}\label{key}
	\widetilde{F}_r(x,z|x_0,x_r)=\widetilde{F}(x,\zeta|x_0)+rz\widetilde{Q}(x,\zeta|x_0)\widetilde{F}_r(x,z|x_r),
\end{equation}
where at the rhs $\widetilde{F}_r(x,z|x_r)$ has the expression given in \eref{eq:FGF_r_x0} by replacing $ x_0 $ with $ x_r $. We thus get:
\begin{equation}\label{eq:FGF_r_xr}
	\widetilde{F}_r(x,z|x_0,x_r)=\frac{(1-z)\widetilde{F}(x,\zeta|x_0)+rz\widetilde{F}(x,\zeta|x_r)}{1-z+rz\widetilde{F}(x,\zeta|x_r)},
\end{equation}
and note that when $ x_0=x_r $ we recover \eref{eq:FGF_r_x0}. In the limit $ z\to1^- $ we can write
\begin{equation}\label{eq:FGF_rxr_exp}
	\widetilde{F}_r(x,z|x_0,x_r)\sim1-(1-z)\left[\frac{1-\widetilde{F}(x,1-r|x_0)}{r\widetilde{F}(x,1-r|x_r)}\right],
\end{equation}
hence we can immediately deduce:
\begin{itemize}
	\item[a)] The process under stochastic resetting is always recurrent independently of the features of the generating function $ \widetilde{F}(x,z|y) $ of the free process, since the leading order term of the expansion \eref{eq:FGF_rxr_exp} is one for any choice of $ x $, $ x_0 $ and $ x_r $;
	\item[b)] In the limit $ z\to1^- $ that we are considering, the generating function in presence of resetting assumes the form dictated by the Lamperti theorem, see \eref{eq:F_form} in section \ref{s:RecTra}, even if in this case we are considering not a first return problem, but the first hitting probability. 
	\item[c)] The term in square brackets at the rhs of \eref{eq:FGF_rxr_exp} is the mean first hitting time $\langle \tau_r(x|x_0,x_r) \rangle$ of the walk with resetting, since it is the coefficient of the linear term in the expansion of $ \widetilde{F}_r(x,z|x_0,x_r) $ in powers of $ 1-z $. Indeed, calling $ F_r(n)=F_r(x,n|x_0,x_r) $ for brevity, we can write
	\begin{equation}\label{key}
		\sum_{n=1}^{\infty}F_r(n)\left[1-(1-z)\right]^n\sim\sum_{n=1}^{\infty}F_r(n)-(1-z)\sum_{n=1}^{\infty}nF_r(n),
	\end{equation}
	in the limit $ z\to1^{-} $. Note that this works because both series appearing in the expansion converge. In principle, one can obtain information on higher order moments of the first hitting time by considering further powers of $ 1-z $. For example the next term would be
	\begin{eqnarray}\label{key}
		\frac{(1-z)^2}2\sum_{n=2}^{\infty}n(n-1)F_r(n)&=&\frac{(1-z)^2}2\left[\mu_2-F_r(1)-\mu_1+F_r(1)\right]\nonumber\\
		&&\frac{(1-z)^2}2\left[\mu_2-\mu_1\right],
	\end{eqnarray}
	where we have used the short notation $ \mu_q=\langle\tau_r^q(x|x_0,x_r)\rangle $. By expanding \eref{eq:FGF_r_xr} up to $ (1-z)^2 $ we get
	\begin{eqnarray}\label{eq:tau2}
		\langle \tau_r^2(x|x_0,x_r)\rangle&=&2\langle\tau_r(x|x_0,x_r)\rangle\Bigg[\frac12-(1-r)\frac{\widetilde{F}'(x,1-r|x_r)}{\widetilde{F}(x,1-r|x_r)}+\nonumber\\
		&&\frac{1-r\widetilde{F}(x,1-r|x_r)}{r\widetilde{F}(x,1-r|x_r)}-\frac{(1-r)\widetilde{F}'(x,1-r|x_0)}{1-\widetilde{F}(x,1-r|x_0)}\Bigg].
	\end{eqnarray}
\end{itemize}

\subsection{Summary of the main results for the Gillis random walk}
Let us summarize the main results for the GRW. As we will explain later in section \ref{s:Results}, depending on the value of $ \epsilon $ defining the jump probability, see \eref{eq:Trans_def}, the walk is transient for $ -1<\epsilon<-\case12 $ and recurrent otherwise. In all cases, we obtain the exact expression of mean first hitting time of the origin $ x=0 $, starting from a generic point $ x_0 $ and with resetting position $ x_r $. We also derive the asymptotics of the mean $ \langle\tau_{r}(0|x_0,x_r)\rangle $ and the standard deviation $ \sigma_r $, in the regimes $ r\to0 $ and $ r\to1 $. The former case is summarized in table \ref{t:t1}, which shows different behaviours depending on the recurrence properties of the walk. In the latter case instead, we find that in the whole range $ -1<\epsilon\leq1 $, hence independently of the recurrence properties, $ \langle\tau_{r}(0|x_0,x_r)\rangle $ diverges as
\begin{equation}\label{key}
	\langle\tau_{r}(0|x_0,x_r)\rangle\sim\frac{|x_r|!\Gamma(1+\epsilon)}{\Gamma(1+\epsilon+|x_r|)}\left(\frac2{1-r}\right)^{|x_r|},\qquad r\to1
\end{equation}
and the standard deviation $ \sigma_r $ follows exactly the same behaviour. By comparing this result with those reported in table \ref{t:t1}, we also conclude that:
\begin{itemize}
	\item[a)] In the transient and null-recurrent regimes, $ \langle\tau_{r}(0|x_0,x_r)\rangle $ diverges for both $ r\to0 $ and $ r\to1 $. This means that resetting is able to expedite the search for any $ 0<r<1 $, and furthermore, there exists an optimal resetting probability $ r^* $ which minimizes the mean first hitting time. In section \ref{s:Results} we consider the particular case $ x_0=x_r $ and present a numerical study on $ r^* $ as a function of $ \epsilon $ and $ x_0 $;
	\item[b)] In the positive-recurrent regime, while $ \langle\tau_{r}(0|x_0,x_r)\rangle $ still diverges for $ r\to1 $, the divergence for $ r\to0 $ disappears, and the mean first hitting time converges instead to a finite positive value, see table \ref{t:t1}. The existence of an optimal resetting probability is thus not guaranteed. Indeed, sticking to $ x_0=x_r $, we show that in the particular case $ \epsilon=1 $ and $ |x_0|=1 $, the resetting only hinders the mean first hitting time. However, if we consider $ \epsilon=1 $ but $ |x_0|>1 $, or $ \epsilon<1 $ with any initial distance, we find that there exists an optimal probability $ r^* $ which minimizes the mean first hitting time. Nevertheless, we explain that, contrarily to the previous ones, in this regime resetting can be advantageous only up to a certain threshold $ r_{\mathrm{th}} $, that we study numerically.
\end{itemize}
Although for the GRW we can obtain exact analytical results, in the following we will also show that the asymptotic behaviours can be derived from a general perspective, that does not require the full knowledge of the underlying generating function $ \widetilde{F}(x,z|x_0) $, but only its general structure, which for a wide class of processes is dictated by the recurrence properties of the system.
\begin{table}[]
	\centering
	\bgroup
	\def\arraystretch{2.2}
	\begin{tabular}{|l|c|c|c|}
		\hline
		& Range of $\epsilon$         & $ \langle\tau_{r}(0|x_0,x_r)\rangle $ ($ r\to0 $) & $ \sigma_r $ ($ r\to0 $)\\ \hline
		Transient                       & $ -1<\epsilon<-\case12 $    & $\displaystyle\frac{C(x_0,x_r)}{r}$            & $ \displaystyle\frac{c(x_0,x_r)}{r} $                                                 \\ \hline
		\multirow{3}{*}{Null-recurrent} & $\epsilon=-\case12$         & $\displaystyle\frac{A(x_0)}{r\log(1/r)}$       & $\displaystyle \sqrt{\frac{A(x_0)}{\log(1/r)}}\frac1r$                                \\ \cline{2-4} 
		& $-\case12<\epsilon<\case12$ & $\displaystyle\frac{B(x_0)}{r^{1/2-\epsilon}}$ & \multirow{3}{*}{$\displaystyle\frac{\sqrt{2\mathcal{B}(x_0)}}{r^{(3/2-\epsilon)/2}}$} \\ \cline{2-3}
		& $\epsilon=\case12$          & $\displaystyle\frac{x_0^2}{2}\log(1/r)$        &                                                                                       \\ \cline{1-3}
		Positive-recurrent              & $\case12<\epsilon\leq1$     & $\displaystyle\frac{x_0^2}{2\epsilon-1}$       &                                                                                       \\ \cline{1-3} \cline{4-4} 
	\end{tabular}
	\egroup
	\caption{Asymptotic (small-$ r $) behaviours for the mean and the standard deviation of the first hitting time, for the GRW with resetting. Here we consider $ x=0 $ with generic starting and resetting positions. The coefficients appearing in the table are defined later in the main text: $ A(x_0) $ is given by \eref{eq:A}, $B(x_0)$ by \eref{B}, $ \mathcal{B}(x_0) $ by \eref{eq:cal_B}, while $ C(x_0,x_r) $ and $ c(x_0,x_r) $ are defined by \eref{eq:C} and \eref{eq:small_c}, respectively.}
	\label{t:t1}
\end{table}

\section{General analysis on the mean first hitting time}\label{s:Anal}
As described by the expressions we found for $ \langle \tau_r(x|x_0,x_r)\rangle $ and $ \langle\tau^2_r(x|x_0,x_r)\rangle $ in \eref{eq:FGF_rxr_exp} and \eref{eq:tau2} respectively, these quantities can be written in terms of the generating functions of the reset-free process. Now, even though the generating functions we are considering here do not refer to the first return probability, it is safe to assume that in the limit $ z\to1^{-} $ they must attain a similar form of that prescribed by Lamperti theorem; namely, for recurrent walks
\begin{equation}\label{eq:F(z|x0)_asym}
	\widetilde{F}(x,z|x_0)\sim1-(1-z)^{\rho}L\left(\frac{1}{1-z};x,x_0\right),
\end{equation}
where in this case the slowly-varying function contains information on the starting point and the target, while we suppose to observe no difference in the value of  $ \rho $.  This is due to the fact that for a random walk where any point $ x $ is accessible from any other point $ x_0 $, the starting point does not change the recurrence properties of the process \cite{Gill-bias}, and we expect that the first return probability and first hitting probability share the same asymptotic behaviour, which is ruled by $ \rho $. More generally, as shown in a previous paper \cite{ROAP}, the exponent $ \rho $ is also related to the temporal scaling of the probability distribution, which can be safely assumed independent of $ x $ and $ x_0 $. This can be more rigorously proved in the context of infinite ergodic theory, see \cite{ROAP,AghKesBar-2020,BarRadAki-2021,Lei-Bar}.

For transient walks, this correspondence between the first return and the first hitting probability is weaker. For example, let us imagine a one-dimensional walk with a strong positive bias, so that the walk drifts towards $ +\infty $: in this case, each point is transient, but the transition from $ x_0 $ to $ x $, if $ x_0<x $, occurs with probability one. Hence, while the first return probability generating function converges to a finite value smaller than one for each choice of the starting point, the generating function of the first hitting probability of $ x $ starting from $ x_0 $ converges instead to one, for $ x_0<x $. On the other hand, for $ x<x_0 $, due to the bias the walk only has probability $ R(x|x_0)<1 $ of reaching $ x $, hence $ \widetilde{F}(x,z|x_0)\to R(x|x_0)<1 $ in the limit $ z\to 1^- $. Since many of the results we derive in the following mainly depend on the whether $ \widetilde{F}(x,z|x_0) $ converges to one or to a smaller value, when we consider transient walks we will only refer to the former case, where the total transition probability from $ x_0 $ to $ x $ occurs with probability $ R(x|x_0)<1 $.

\subsection{High resetting probability}
Let us  first examine
\begin{equation}\label{eq:tau0r_disp}
	\langle\tau_r(x|x_0,x_r)\rangle=\frac{1-\widetilde{F}(x,1-r|x_0)}{r\widetilde{F}(x,1-r|x_r)},
\end{equation}
in the regime $ r\approx1 $, so that we have to consider the behaviour of $ \widetilde{F}(x,z|y) $ close to $ z=0 $. We recall that we are only considering walks that can perform nearest-neighbour jumps, hence it is easy to see that the generating function must vanish as $ z^{|x-y|} $ in the limit. This can be seen by observing that the generating function is
\begin{equation}\label{key}
	\widetilde{F}(x,z|y)=\sum_{n=1}^{\infty}F(x,n|y)z^n,		
\end{equation}
and that due to the constraint of nearest-neighbour jumps, in general the particle can not cover a distance of size $ d(x,y)=|x-y| $ in less than $ d(x,y) $ steps. Therefore, for $ m=1,2,\dots,d(x,y)-1$, all coefficients $ F(x,m|y) $ must vanish, and we can write
\begin{equation}\label{key}
	\widetilde{F}(x,z|y)=z^{d(x,y)}\sum_{m=0}^{\infty}F(x,d(x,y)+m|y)z^m,	
\end{equation}
where the simple shift of index $ m=n-d(x,y) $ has been considered. This shows that as $ z\to0 $ the function $ \widetilde{F}(x,z|y) $ is asymptotic to
\begin{equation}\label{key}
	 \widetilde{F}(x,z|y) \sim \frac{z^{|x-y|}}{K(x,y)},
\end{equation}
thus yielding
\begin{equation}\label{eq:tau_asym_r1}
	\langle\tau_r(x|x_0,x_r)\rangle\sim K(x,x_r)(1-r)^{-|x-x_r|},
\end{equation}
where the coefficient $ K(x,x_r) $ has been introduced just to simplify the notation and is precisely given by
\begin{equation}\label{key}
	K(x,x_r)=\frac{1}{F(x,d(x,x_r)|x_r)}.
\end{equation}
In particular, for the simple and symmetric random walk, it is easy to see that the probability of reaching $ x $ from $ x_r $ for the first time in exactly $ |x-x_r| $ steps is given by $ F(x,d(x,x_r)|x_r)=2^{-d(x,x_r)} $, hence one obtains
\begin{equation}\label{key}
	\langle\tau_r(x|x_0,x_r)\rangle\sim \left(\frac{2}{1-r}\right)^{|x-x_r|}.
\end{equation}
We can immediately see that in this regime the mean first hitting time is not affected by the starting point. This is obviously a trivial observation, because in the regime of high reset probability, the information on $ x_0 $ is immediately lost. Also, we observe that this limiting behaviour does not change if we consider a recurrent or a transient process. Note that here we have a power-law divergence which depends only on the distance between the resetting position and $ x $.

\subsection{Low resetting probability}
We now consider $ r\approx0 $, so that the behaviour of $ \widetilde{F}(x,z|y) $ close to $ z=1 $ is important. Let us begin with the recurrent case, where we assume \eref{eq:F(z|x0)_asym} to be valid, hence by substituting it in \eref{eq:tau0r_disp} we get
\begin{equation}\label{eq:tau_asym_r0_rec}
	\langle\tau_r(x|x_0,x_r)\rangle\sim\frac{r^\rho L(1/r;x,x_0)}{r-r^{1+\rho}L(1/r;x,x_r)}\sim\frac1{r^{1-\rho}}L\left(\frac1r;x,x_0\right),	
\end{equation}
We see that in this regime only $ x_0 $ matters, as expected. As long as $ \rho<1 $, we have a power-law divergence which depends only on the characteristics of the reset-free process, i.e. by the exponent $ \rho $. The combination of two divergences at both ends of the range of the resetting parameter is the typical behaviour observed for the mean first passage time in the presence of resetting, for both diffusion and random walks, see \cite{EvaMaj-2011,Reu-2016,PalKunEva-2016,PalReu-2017,CheSok-2018,TalPalSek-2020,BesBovPet-2020,SinPal-2022,RAD-2022,BonPal-2021,RiaBoyHer-2020,RiaBoyHer-2022}. As a consequence, one finds that the mean first passage time curve has a minimum in correspondence of the optimal reset value\textemdash note that for diffusion processes, one considers the reset rate rather then the resetting probability, which takes values in $ (0,\infty) $. For $ \rho=1 $ instead, the mean first hitting time may either diverge as a slowly-varying function, or converge to a constant, which must coincide with the value of the mean first hitting time in absence of resetting. In this case, it is not trivial to deduce the existence of an optimal resetting value, and more generally, it is not guaranteed that resetting can expedite the search process \textit{for any value} of the resetting parameter, as instead happens for $ \rho<1 $ (remarkably, this observation is the basis of many recent works \cite{BonPal-2021,RiaBoyHer-2022,RayMonReu-2019,RayReu-2020,RayReu-2021,MerBoyMaj-2022,Ray-2022}).

We can now repeat the analysis for the transient case. As pointed out in section \ref{s:RecTra}, we consider the case where the generating function $ \widetilde{F}(x,z|y) $ converges to a constant in the limit $ z\to1^{-} $, which now depends on both $ x $ and $ y $:
\begin{equation}\label{key}
	\lim_{z\to1^{-}}\widetilde{F}(x,z|y)=R(x|y).
\end{equation}
Therefore, as $ r\to0 $, from \eref{eq:tau0r_disp} we have
\begin{equation}\label{eq:tau_asym_r0_tr}
	\langle\tau_r(x|x_0,x_r)\rangle\sim\frac{Q(x|x_0)}{R(x|x_r)}\frac1r,
\end{equation}
where we recall that $ Q(x|y)=1-R(x|y) $ is the probability of never hitting $ x $ starting from $ y $. We observe that, contrarily to the previous regime, in the transient case the information on the resetting position is not lost, even in the limit $ r\to0 $. This happens because the reset-free process has a positive probability of never hitting the target, hence the walks that are reset at least once remain statistically relevant even in the limit of very low resetting probability.

To summarize, the asymptotic behaviour of the mean first hitting time for small values of the resetting probability, viz. $ r\approx0 $, is given by
\begin{equation}\label{key}
	\langle\tau_r(x|x_0,x_r)\rangle\sim
	\cases{\frac{Q(x|x_0)}{R(x|x_r)}\,\frac1r & transient walk\\
	\frac{L(1/r;x,x_0)}{r^{1-\rho}}& recurrent walk $(0\leq\rho\leq1)$.
	}
\end{equation}
In particular, for $ \rho=\case12 $, which corresponds to the simple and symmetric random walk, one recovers the behaviour $ 1/\sqrt{r} $ for small $ r $, which is also observed for diffusion processes \cite{EvaMaj-2011}.

\section{Second moment and standard deviation}\label{s:Anal2}
From the results of the previous subsection, we can easily deduce the behaviour of the second moment for $ r\approx1 $ and $ r\approx0 $. It is convenient to consider the ratio
\begin{eqnarray}\label{eq:tau2/tau}
	\frac{\langle \tau_r^2(x|x_0,x_r)\rangle}{2\langle\tau_r(x|x_0,x_r)\rangle}&=&\frac{1-r\widetilde{F}(x,1-r|x_r)}{r\widetilde{F}(x,1-r|x_r)}-(1-r)\frac{\widetilde{F}'(x,1-r|x_r)}{\widetilde{F}(x,1-r|x_r)}+\nonumber\\
	&&\frac12-(1-r)\frac{\widetilde{F}'(x,1-r|x_0)}{1-\widetilde{F}(x,1-r|x_0)}.
\end{eqnarray}

\subsection{High resetting probability}
Once again, let us consider first $ r\approx1 $. Since
\begin{eqnarray}\label{key}
	\widetilde{F}(x,z|y)&\sim&\frac{z^{|x-y|}}{K(x,y)}\\
	\widetilde{F}'(x,z|y)&\sim&\frac{|x-y|}{K(x,y)}z^{|x-y|-1},
\end{eqnarray}
the leading-order term is the first one at the rhs of \eref{eq:tau2/tau}, hence
\begin{equation}\label{key}
	\frac{\langle \tau_r^2(x|x_0,x_r)\rangle}{2\langle\tau_r(x|x_0,x_r)\rangle}\sim\frac1{\widetilde{F}(x,1-r|x_r)}\sim K(x,x_r)(1-r)^{-|x-x_r|},
\end{equation}
which means $ \langle\tau_r^2(x|x_0,x_r)\rangle\sim2\langle\tau_r(x|x_0,x_r)\rangle^2 $, see \eref{eq:tau_asym_r1}. If we introduce the standard deviation
\begin{equation}\label{eq:std_dev_disp}
	\sigma_r=\sqrt{\langle\tau_r^2(x|x_0,x_r)\rangle-\langle\tau_r(x|x_0,x_r)\rangle^2},
\end{equation}
then we obtain
\begin{equation}\label{eq:sig_asym_r1}
	\sigma_r\sim\langle\tau_r(x|x_0,x_r)\rangle\sim K(x,x_r)(1-r)^{-|x-x_r|},
\end{equation}
which indicates that the mean and the standard deviation of the mean first hitting time share the same divergence in the limit $ r\to1 $, suggesting that the distribution of $ \tau_{r}(x|x_0,x_r) $ becomes approximately exponential.

\subsection{Low resetting probability}
In the opposite limit, viz. when $ r\approx0 $, we have to evaluate the behaviour of the derivative $ \widetilde{F}'(x,z|y) $. Let us begin with the recurrent case, where we recall
\begin{equation}\label{eq:F_Lamp_II}
	\widetilde{F}(x,z|y)\sim1-(1-z)^\rho L\left(\frac1{1-z};x,y\right),
\end{equation}
with $ 0\leq\rho\leq1$. We have to consider the behaviour of the derivative $ \widetilde{F}'(x,z|y) $ close to $ z=1 $, and we can say
\begin{equation}\label{eq:F_der}
	\widetilde{F}'(x,z|y)\sim\frac{\rho}{(1-z)^{1-\rho}} L\left(\frac1{1-z};x,y\right),
\end{equation}
see the original paper by Lamperti \cite{Lam}. Hence we obtain from \eref{eq:tau2/tau} that the leading order is given by
\begin{equation}\label{eq:LD_sig}
	\frac{\langle \tau_r^2(x|x_0,x_r)\rangle}{2\langle\tau_r(x|x_0,x_r)\rangle}\sim\frac{1-\rho}r,
\end{equation}
thus $ \langle\tau_r^2(x|x_0,x_r)\rangle\sim2(1-\rho)\langle\tau_r(x|x_0,x_r)\rangle/r $. By considering \eref{eq:tau_asym_r0_rec}, we obtain that the standard deviation diverges as
\begin{equation}\label{eq:sig_asym_r0_rec}
	\sigma_r\sim\sqrt{\frac{2}{r}(1-\rho)\langle\tau_r(x|x_0,x_r)\rangle}\sim\sqrt{2(1-\rho)L(1/r;x,x_0)}\,r^{-(1-\rho/2)}.
\end{equation}
Note that for $ \rho=1 $ we need to be particularly cautious, because the leading-order term that we have computed in \eref{eq:LD_sig} vanishes, and we have thus to consider the first subleading correction. In other words, we have to evaluate the subleading terms in the asymptotic expansion of the slowly-varying function. As a first approximation, by taking the derivative with respect to $ z $ of \eref{eq:F_Lamp_II}, with $ \rho=1 $, we can say
\begin{equation}\label{key}
	\widetilde{F}'(x,z|y)\sim L\left(\frac1{1-z};x,y\right)-\frac1{1-z}L'\left(\frac1{1-z};x,y\right),
\end{equation}
and the first correction may arise from both terms at the rhs. We report the details in \ref{a:Fr_der}, here we only state the main results. We distinguish three cases, depending on the statistical features of the underlying first hitting process:
\begin{enumerate}
\item \textit{The mean diverges as the slowly-varying function} $ L(t;x,y) $: in this scenario, for $ t\to\infty $ we can write
\begin{equation}\label{key}
	L'(t;x,y)\sim\frac{\ell(t;x,y)}{t},
\end{equation}
where $ \ell(t) $ is a slowly-varying function. Note that $ \ell(t;x,y) $ must be subdominant with respect to $ L(t;x,y) $, see also \ref{a:Fr_der}. Then from \eref{eq:tau2/tau} we get
\begin{equation}\label{key}
	\frac{\langle\tau_r^2(x|x_0,x_r)\rangle}{2\langle\tau_r(x|x_0,x_r)\rangle}\sim\frac{\ell(1/r;x,x_0)}{rL(1/r;x,x_0)},
\end{equation}
and by recalling that in this case $\langle\tau_r(x|x_0,x_r)\rangle\sim L(1/r;x,x_0)$, we then have
\begin{equation}\label{key}
	\sigma_r\sim\sqrt{2\ell(1/r;x,x_0)}r^{-1/2}.
\end{equation}
\item \textit{The mean is finite, but the second moment is infinite}: in this situation the first hitting time probability decays as $ F(x,n|y)\propto n^{-2-\delta} $, with $ 0<\delta<1 $, and we must have
\begin{equation}\label{key}
	L(t;x,y)\to\mu_1(x,y)\qquad\textrm{as }t\to\infty,
\end{equation}
where $ \mu_1(x,y) $ is the mean first hitting time of the underlying process. In \ref{a:Fr_der} we find that 
\begin{eqnarray}
	L(t;x,y)&\sim&\mu_1(x,y)-\frac{1}{t^{\delta}}\ell(t;x,y)\\
	L'(t;x,y)&\sim&\frac{\delta}{t^{1+\delta}}\ell(t;x,y)\label{eq:L'_ii},
\end{eqnarray}
thus from \eref{eq:tau2/tau} we obtain
\begin{equation}\label{key}
	\frac{\langle\tau_r^2(x|x_0,x_r)\rangle}{2\langle\tau_r(x|x_0,x_r)\rangle}\sim\frac{\ell(1/r;x,x_0)}{\mu_1(x,x_0)}\frac\delta{r^{1-\delta}},
\end{equation}
and by recalling now that $\langle\tau_r(x|x_0,x_r)\rangle\to \mu_1(x,x_0)$, we arrive at
\begin{equation}\label{key}
	\sigma_r\sim\sqrt{2\delta\ell(1/r;x,x_0)}r^{-(1-\delta)/2}.
\end{equation}
\item \textit{Both the mean and the second moment are finite}: let us call $ \mu_1(x,y) $ the first moment of the first hitting time distribution of $ x $ starting from $ y $ for the underlying process, and $ \mu_2(x,y) $ the second moment. In this case we find
\begin{eqnarray}\label{key}
	L(t;x,y)&\sim&\mu_1(x,y)-\left[\mu_2(x,y)-\mu_1(x,y)\right]\frac1{2t}\\
	L'(t;x,y)&\sim&\left[\mu_2(x,y)-\mu_1(x,y)\right]\frac1{2t^{2}},
\end{eqnarray}
then from \eref{eq:tau2/tau}
\begin{equation}\label{key}
	\frac{\langle\tau_r^2(x|x_0,x_r)\rangle}{2\langle\tau_r(x|x_0,x_r)\rangle}\sim\frac{\mu_2(x,x_0)}{2\mu_1(x,x_0)},
\end{equation}
where, as in the previous case, $\langle\tau_r(x|x_0,x_r)\rangle\to \mu_1(x,x_0)$. Therefore the standard deviation in presence of resetting converges to that of the reset-free process
\begin{equation}\label{key}
	\sigma_r\to\sigma_0.
\end{equation}
\end{enumerate}

We can now tackle the transient case, with $ \widetilde{F}(x,z|y)\sim R(x|y) $ as $ z\to1^- $. It is easy to see that this time we get
\begin{equation}
	\frac{\langle \tau_r^2(x|x_0,x_r)\rangle}{2\langle\tau_r(x|x_0,x_r)\rangle}\sim\frac{1}{rR(x|x_r)},
\end{equation}
and thus we obtain that the standard deviation behaves as
\begin{equation}\label{eq:sig_asym_r0_tr}
	\sigma_r\sim\langle\tau_r(x|x_0,x_r)\rangle\sqrt{\frac{2}{Q(x|x_0)}-1}\sim\sqrt{\frac{1-R^2(x|x_0)}{R^2(x|x_r)}}\,\frac1r.
\end{equation}
It follows that in the limit $ r\to0 $, the ratio $ \sigma_r/\langle\tau_r(x|x_0,x_r)\rangle $ is always a constant bigger than one, indeed
\begin{equation}\label{key}
	\frac{\sigma_r}{\langle\tau_r(x|x_0,x_r)\rangle}\sim\sqrt{\frac{1-R^2(x|x_0)}{R^2(x|x_r)}}\cdot\frac{R(x|x_r)}{1-R(x|x_0)}=\sqrt{\frac{1+R(x|x_0)}{1-R(x|x_0)}}.
\end{equation}

To summarize, we found that for small values of the resetting probability, viz. $ r\approx0 $, we can write
\begin{equation}\label{key}
	\sigma_r\sim
	\cases{\sqrt{\frac{1-R^2(x|x_0)}{R^2(x|x_r)}}\,\frac1r & transient walk\\
		\sqrt{2(1-\rho)L\left(\frac1r;x,x_0\right)}\,\frac1{r^{1-\rho/2}} & recurrent walk $(0\leq\rho<1)$,
	}
\end{equation}
and in particular $ \sigma_r\propto r^{-3/4} $ for $ \rho=\case12 $, while for $ \rho=1 $, the value of another exponent $ \delta$ is necessary to determine the behaviour of $ \sigma_r $. This exponent is related to the decay of the first hitting probability of the underlying process, which we indicate as $ F(x,n|y)\propto n^{-2-\delta} $, and we obtain
\begin{equation}\label{key}
	\sigma_r\sim
	\cases{\sqrt{2\ell(1/r;x,x_0)}\,\frac1{r^{1/2}} & for $\delta=0$\\
		\sqrt{2\delta\ell(1/r;x,x_0)}\,\frac1{r^{(1-\delta)/2}} & for $0<\delta\leq1$,
	}
\end{equation}
while for $ \delta>1 $ the standard deviation converges to that of the underlying first passage process, which indeed in this case is finite, and so we have
\begin{equation}\label{key}
	\sigma_r\to\sigma_0.
\end{equation}

\section{Optimal resetting}\label{s:Opt}
According to what we have seen so far, we can relate the recurrence properties of the underlying process, in particular, the value of the exponent $ \rho $, to the effects of resetting on the mean first hitting time. From the previous discussion, it is evident that for both transient and null-recurrent walks the introduction of resetting expedites the first hitting process \textit{for any value} of the resetting probability $ 0<r<1 $. Indeed, as already thoroughly discussed in the literature, in these cases the mean first hitting time curve presents a divergence at both ends of this range. While the divergence close to $ r\approx1 $ is intrinsically related to the resetting mechanism, the one for $ r\approx0 $ depends on the behaviour of the reset-free system; hence, for transient and null-recurrent walks, the introduction of a finite mean first hitting time for any $ 0<r<1 $ makes resetting always advantageous. Furthermore, the non-monotonic behaviour of the curve assures that there is an optimal value of the resetting probability $ r^* $ which minimizes the mean first hitting time. This is obviously the solution of
\begin{equation}\label{eq:der=0}
	\frac{\mathrm{d}}{\mathrm{d}r^{*}}\langle\tau_{r^*}(x|x_0,x_r)\rangle=0.
\end{equation}

In general, it has been pointed out \cite{BonPal-2021,RiaBoyHer-2022} that to evaluate whether resetting may or not be beneficial it is useful to introduce the \textit{coefficient of variation}, which is defined in terms of the first and second moment of the first hitting time distribution
\begin{equation}\label{key}
	\xi_r(x|x_0,x_r)=\frac{\sqrt{\langle\tau^2_r(x|x_0,x_r)\rangle-\langle\tau_r(x|x_0,x_r)\rangle^2}}{\langle\tau_r(x|x_0,x_r)\rangle},
\end{equation}
and for which a more explicit expression can be provided by considering \eref{eq:tau2} and \eref{eq:tau0r_disp}. Let us first observe that, by defining
\begin{equation}\label{key}
	\mathcal{F}(x,r|x_0,x_r)=\frac{1-\widetilde{F}(x,1-r|x_0)}{\widetilde{F}(x,1-r|x_r)},
\end{equation}
we have
\begin{eqnarray}
	\frac{\mathrm{d}}{\mathrm{d}r}\log\bigg(\mathcal{F}(x,r|x_0,x_r)\bigg)&=&\frac{\mathrm{d}}{\mathrm{d}r}\log\bigg(\langle\tau_{r}(x|x_0,x_r)\rangle\bigg)+\frac1r\label{eq:logF_cal}\\
	&=&\frac{\widetilde{F}'(x,1-r|x_0)}{1-\widetilde{F}(x,1-r|x_0)}+\frac{\widetilde{F}'(x,1-r|x_r)}{\widetilde{F}(x,1-r|x_r)}.
\end{eqnarray}
This is useful to rewrite the second moment, since the rhs of the equation above explicitly appears in \eref{eq:tau2}, and with some straightforward algebra we then obtain
\begin{equation}\label{eq:tau2_s7}
	\fl\langle\tau_{r}^2(x|x_0,x_r)\rangle=2\langle\tau_{r}(x|x_0,x_r)\rangle\left[\langle\tau_{r}(x|x_r)\rangle+\frac12\right]-2(1-r)\frac{\mathrm{d}}{\mathrm{d}r}\langle\tau_{r}(x|x_0,x_r)\rangle,
\end{equation}
where
\begin{equation}\label{key}
	\langle\tau_{r}(x|x_r)\rangle=\langle\tau_{r}(x|x_r,x_r)\rangle=\frac{1-\widetilde{F}(x,1-r|x_r)}{r\widetilde{F}(x,1-r|x_r)}
\end{equation}
is the mean first hitting time for $ x_r=x_0 $. It follows that the square coefficient of variation satisfies for all $ 0<r<1 $ the equation
\begin{equation}\label{eq:zeta}
	\xi_{r}^2(x|x_0,x_r)+1=\frac{2\langle\tau_{r}(x|x_r)\rangle+1}{\langle\tau_{r}(x|x_0,x_r)\rangle}+2(1-r)\frac{\mathrm{d}}{\mathrm{d}r}\left(\frac1{\langle\tau_{r}(x|x_0,x_r)\rangle}\right).
\end{equation}
It has been argued that a sufficient condition for resetting to be advantageous is that the derivative of $ \langle\tau_r(x|x_0,x_r)\rangle $ (with respect to $ r $) is negative for $ r $ sufficiently small \cite{BonPal-2021,RiaBoyHer-2022}, and this can be translated to a simple inequality involving $ \xi_{r}(x|x_0,x_r) $. To see this, just note that if for some $ r $ we have
\begin{equation}\label{eq:der_neg}
	\frac{\mathrm{d}}{\mathrm{d}r}\langle\tau_r(x|x_0,x_r)\rangle <0,
\end{equation}
then it follows from \eref{eq:zeta}:
\begin{equation}\label{eq:zeta_ineq}
	\xi_{r}^2(x|x_0,x_r)+1>\frac{2\langle\tau_{r}(x|x_r)\rangle+1}{\langle\tau_{r}(x|x_0,x_r)\rangle}.
\end{equation}
It is then sufficient to check if this inequality is valid for $ r $ small enough. Furthermore, note that whenever the derivative of $\langle\tau_{r}(x|x_0,x_r)\rangle$ is positive, the opposite inequality holds. Then when the derivative vanishes, viz. for $ r=r^* $, we must have
\begin{equation}\label{key}
	\xi_{r^*}^2(x|x_0,x_r)+1=\frac{2\langle\tau_{r^*}(x|x_r)\rangle+1}{\langle\tau_{r^*}(x|x_0,x_r)\rangle},
\end{equation}
which establishes a condition for the optimal value $ r^* $. We remark that for $ x_r=x_0 $ one recovers the same equation found by Riascos \textit{et al}. \cite{RiaBoyHer-2022}:
\begin{equation}\label{eq:CV_opt0}
	\xi_{r^*}^2(x|x_0)=\frac{1}{\langle\tau_{r^*}(x|x_0)\rangle}+1.
\end{equation}

As we have already said, for transient and null-recurrent walks, it is unnecessary to check if \eref{eq:zeta_ineq} is verified, since in both cases the mean first hitting time diverges for $ r\to0 $. However, the situation may change significantly if we consider positive-recurrent walks, since in this case for $ r\approx0 $ the mean first hitting time reaches a finite value. Let us first consider the case of positive-recurrent walks with infinite variance; then we recall
\begin{equation}\label{key}
	\widetilde{F}(x,z|y)\sim1-(1-z)\mu_1(x,y)+(1-z)^{1+\delta}\ell\left(\frac{1}{1-z};x,y\right),
\end{equation}
hence
\begin{equation}\label{key}
	\langle\tau_r(x|x_0,x_r)\rangle\sim\mu_1(x,x_0)-(1+\delta)\ell\left(\frac1r;x,x_0\right)r^\delta,
\end{equation}
where, from the discussion of the previous section, $ \ell(1/r;x,x_0) $ must be positive as $ r\to0 $. It follows that in this case the mean first hitting time has negative slope for $ r $ small enough, hence resetting is still able to expedite the process. However, differently from the previous case, this does not happen \textit{for any} value of the resetting probability: since $ \langle\tau_r(x|x_0,x_r)\rangle $ diverges for $ r\to1 $, there must be a range of $ r  $ such that
\begin{equation}\label{key}
	\langle\tau_r(x|x_0,x_r)\rangle>\mu_1(x,x_0),
\end{equation}
meaning that the resetting-free process possesses a lower mean first hitting time. Hence in this case we must be aware of the presence of a threshold $ r_{\mathrm{th}} $ above which the resetting increases the mean completion time, which is obviously the positive solution of
\begin{equation}\label{key}
	\langle\tau_{r_{\mathrm{th}}}(x|x_0,x_r)\rangle=\mu_1(x,x_0).
\end{equation}

Finally, let us consider the case of underlying first hitting processes with also finite variance, where we can write
\begin{equation}\label{key}
	\widetilde{F}(x,z|y)\sim1-(1-z)\mu_1(x,y)+\frac12(1-z)^2\left[\mu_2(x,y)-\mu_1(x,y)\right],
\end{equation}
which yields
\begin{equation}\label{key}
	\langle\tau_{r}(x|x_0,x_r)\rangle\sim\mu_1(x,x_0)-\frac r2\left[\mu_2(x,x_0)-\mu_1(x,x_0)-2\mu_1(x,x_0)\mu_1(x,x_r)\right].
\end{equation}
Therefore the mean first hitting time has negative slope for $ r\approx0 $ if the moments of the reset-free process satisfy
\begin{equation}\label{key}
	\mu_2(x,x_0)>2\mu_1(x,x_0)\mu_1(x,x_r)+\mu_1(x,x_0),
\end{equation}
or
\begin{equation}\label{eq:CV_fvar}
	\xi_0^2(x|x_0)+1>\frac{2\mu_1(x,x_r)+1}{\mu_1(x,x_0)},
\end{equation}
where
\begin{equation}\label{key}
	\xi_0(x|x_0)=\frac{\sqrt{\mu_2(x,x_0)-\mu_1^2(x,x_0)}}{\mu_1(x,x_0)}
\end{equation}
is the coefficient of variation of the reset-free process. As remarked in previous works \cite{BonPal-2021}, in this case the criterion only involves the statistical properties of the underlying process; in particular, just the first and the second moment.

We remark that even if the previous condition is not satisfied, the existence of a range of $ r $ for which resetting becomes advantageous can not be excluded, nor the presence of a global minimum in correspondence of an optimal value $ r^* $, see \cite{BonPal-2021,FlyPil-2021} for a few illustrative examples. In other words, the mean first hitting time is not necessarily a monotonic function of $ r $, as it can show for example a maximum and then a minimum below $ \mu_1(x,x_0) $; in other words, the situation becomes much more complex and specific analysis for each case ought to be considered.

\section{Results for the Gillis random walk and numerical simulations}\label{s:Results}
We now apply the results obtained in the previous sections to the GRW and test our findings with numerical simulations. To offer a complete overview of the problem, we first briefly review the properties of the system without resetting, which will be used later to derive the quantities of interest under geometric restart.
\subsection{GRW without resetting}
Let us begin with the return probability. Remarkably, for the GRW the generating function $ \widetilde{F}(z) $ can be computed exactly \cite{Gill-bias}:
\begin{equation}\label{eq:FGF_Ret}
	\widetilde{F}(z)=1-\frac{\hyp\left(\frac12\epsilon,\frac12+\frac12\epsilon;1;z^2\right)}{\hyp\left(\frac12+\frac12\epsilon,1+\frac12\epsilon;1;z^2\right)},
\end{equation}
where $ \hyp\left(a,b;c;t\right) $ is the Gauss hypergeometric function, defined by
\begin{equation}\label{key}
	\hyp\left(a,b;c;t\right)=\sum_{n=0}^{\infty}\frac{(a)_n(b)_n}{(c)_n}\frac{t^n}{n!}.
\end{equation}
Here $ (\nu)_n $ is the Pochhammer symbol, which can be defined in terms of the Gamma function as 	$(\nu)_n=\Gamma(\nu+n)/\Gamma(\nu)$. By taking the limit $ z\to1^{-} $ of \eref{eq:FGF_Ret}, one can show that the walk is transient for $ -1<\epsilon<-\frac12 $ and recurrent for $ \epsilon\geq-\frac12 $. When the walk is recurrent, from \eref{eq:FGF_Ret} we can get the value of $ \rho $ \cite{OPRA,ROAP}:
\begin{equation}\label{eq:rho_Gillis}
	\rho=\cases{0 & for $\epsilon=-\frac12$\\
		\frac12+\epsilon & for $-\frac12<\epsilon<\frac12$\\
		1 & for $\frac12\leq\epsilon<1$.}
\end{equation}
It is also possible to obtain the large-$ n $ behaviour of the mean return time $ \langle\tau\rangle $ for any value of $ \epsilon>-\frac12 $, which reads
\begin{equation}\label{eq:tau_free}
	\langle\tau\rangle\propto
	\cases{
		\frac{n}{\log^2(n)} &for $\epsilon=-\frac12$\\
		n^{1/2-\epsilon}&for $-\frac12<\epsilon<\frac12$\\
		\log (n)&for $\epsilon=\frac12$,}
\end{equation}
where $ \tau\propto g(n) $ indicates that $ \lim_{n\to\infty}\langle\tau\rangle/g(n) $ is a constant, while in the positive-recurrent case $ \langle\tau\rangle $ converges to the expected value
\begin{equation}\label{key}
	\langle\tau\rangle\to\frac{2\epsilon}{2\epsilon-1},\qquad\epsilon>\frac12.
\end{equation}
In the transient case instead we can explicitly compute the return probability $ R $, given by
\begin{equation}\label{key}
	R=-1-\frac1\epsilon,\qquad\epsilon<-\frac12.
\end{equation}

Let us mention here that recurrence and transience can be explained in physical terms by the presence of an external field, which either forces the particle to repeatedly return to the starting point or eventually drift away from it. Indeed, as we mentioned in the introduction and shown in previous work \cite{OPRA}, the GRW can be interpreted as the microscopic description of a Brownian particle diffusing in a logarithmic potential $ V(x)\sim V_0\log(|x|) $. Let us now drop the assumption of unitary lattice spacing and time step length, and introduce a microscopic jump size $ \ell $ and step duration $ \delta $. We call $ P_{n\delta}(x\ell|x_0\ell) $ the probability of occupying $ z=x\ell $ at time $ t=n\delta $, having started from $ z_0=x_0\ell $. By recalling the definition of the transition probability $ G_{x,y} $ given in \eref{eq:Trans_def}, the evolution equation can be written
\begin{eqnarray}
	\lefteqn{\frac{P_{(n+1)\delta}(x\ell|x_0\ell)-P_{n\delta}(x\ell|x_0\ell)}{\delta} = } \nonumber\\
	& & \frac{\ell^2}{2\delta}\Bigg[\frac{P_{n\delta}((x+1)\ell|x_0\ell)-2P_{n\delta}(x\ell|x_0\ell) +P_{n\delta}((x-1)\ell|x_0\ell) }{\ell^2}\nonumber\\
	&&+\frac{\epsilon P((x+1)\ell|x_0\ell)}{(x+1)\ell^2}-\frac{\epsilon P((x-1)\ell|x_0\ell)}{(x-1)\ell^2}\Bigg],
\end{eqnarray}
which is a finite difference approximation of a partial differential equation. Indeed, the first line can be approximated by a derivative with respect to the time variable $ t=n\delta$; the second line, by a second derivative with respect to the variable $ z=x\ell $; the third line by the derivative of $ \epsilon P(z|z_0)/z $ with respect to $ z $ (up to a factor $ 2 $), see \cite{Olv}. Therefore, by considering the limit $ \ell\to0 $, $ \delta\to0 $, with the ratio $ \ell^2/(2\delta)=D $ kept fixed, we formally get the Fokker-Planck equation:
\begin{equation}\label{eq:Gillis_FP}
	\frac{\partial P(z,t|z_0)}{\partial t}=D\frac{\partial^2 P(z,t|z_0)}{\partial z^2}+2\epsilon D\frac{\partial}{\partial z}\left[\frac{P(z,t|z_0)}{z}\right].
\end{equation}
We recall that the Fokker-Planck equation of a one-dimensional diffusive particle, with diffusion constant $ D $, in a potential $ V(z) $ is written
\begin{equation}\label{key}
	\frac{\partial P(z,t|z_0)}{\partial t}=D\frac{\partial^2 P(z,t|z_0)}{\partial z^2}+\frac{\partial}{\partial z}\left[P(z,t|z_0)\frac{\rmd V(z)}{\rmd z}\right],
\end{equation}
which means that we can interpret \eref{eq:Gillis_FP} with $ V(z) $ of the form $ V(z)=2\epsilon D\log(|z|/a) $, where $ a $ is a cut-off which can be tuned to regularize the potential around the origin. This reveals that the parameter $ \epsilon $ is equal to the ratio $ \epsilon=V_0/(2D) $. It is then intuitive that, by tuning the sign and the relative strength of the potential with respect to the diffusion constant, one can drastically change the recurrence properties of the system: for $ \epsilon>\case12 $, i.e. $ V_0>D $, the potential is sufficiently confining and the process is thus positive-recurrent; in the opposite case, $ V_0<-D $, the particle is drifted away from the origin strongly enough that the process is transient; in the intermediate case, $ -D<V_0<D $, diffusion is instead in control and the walk is hence null-recurrent.

Under these premises, we now turn to the quantity of interest for the present work, which is the first hitting time probability of site $ x $ starting the walk from $ x_0 $. Unfortunately this problem has not been solved yet for any possible $ x $, except for $ x=0 $, hence from now on we will consider only this case. By denoting with $ \widetilde{F}(0,z|x_0) $ the corresponding generating function, we have for any $ x_0\neq0 $
\begin{equation}\label{eq:FGF}
	\widetilde{F}(0,z|x_0)=\frac{z^{|x_0|}}{|x_0|!}\frac{\Gamma(1+\epsilon+|x_0|)}{2^{|x_0|}\Gamma(1+\epsilon)}\,\frac{\hyp\left(\frac {1+\epsilon+|x_0|}2,1+\frac {\epsilon+|x_0|}2;1+|x_0|;z^2\right)}{\hyp\left(\frac {1+\epsilon}2,1+\frac 12\epsilon;1;z^2\right)},
\end{equation}
where $ \Gamma(x) $ is the Gamma function. The interested reader may find a derivation of this result in \ref{a:F_app}. For particular values of $ \epsilon $, the expression above can be written in terms of elementary functions: for $ \epsilon=0 $, the model becomes the simple symmetric random walk on $ \mathbb{Z} $ and indeed one recovers the well-known formula for the generating function of the first hitting time probability; more interestingly, for $ \epsilon=1 $ one can write
\begin{equation}\label{eq:FGF_ep1}
	\widetilde{F}(0,z|x_0)=\left(\frac{z}{1+\sqrt{1-z^2}}\right)^{|x_0|}\left(1+|x_0|\sqrt{1-z^2}\right),
\end{equation}
as originally found in \cite{ChaHug-1988}, where this case was related to the problem of ion diffusion in a semi-infinite domain near a uniformly charged surface; for $ \epsilon=-1 $ the rhs in \eref{eq:FGF} vanishes, and indeed this is not a case we will deal with, because the first hitting time of the origin is pathologically always infinite: it is clear from the definition of the transition probability \eref{eq:Trans_def} that whenever the particle reaches $ x=\pm1 $, it is then pushed away from the origin with probability one.

Of course, one can check that the assumption about the limit form of $ \widetilde{F}(0,z|x_0) $ for $ z\approx1 $ we made in the previous section holds for any $ \epsilon\geq-\case12 $, and verify that the value of $ \rho $ is given again by \eref{eq:rho_Gillis}, see \ref{app:SV}. From this generating function we can compute the expected first hitting time of the origin $ \langle\tau(0|x_0)\rangle $, and as we have already mentioned, we expect a similar behaviour as the mean return time $ \langle\tau\rangle $. Indeed we get
\begin{equation}\label{eq:tau_x0_free}
	\langle\tau(0|x_0)\rangle\sim
	\cases{
			\frac{A(x_0)n}{\log^2(n)} &for $\epsilon=-\frac12$\\
			\frac{(1/2+\epsilon)B(x_0)}{\Gamma(3/2-\epsilon)}n^{1/2-\epsilon}&for $-\frac12<\epsilon<\frac12$\\
			\frac{x_0^2}{2}\log (n)&for $\epsilon=\frac12$,}
\end{equation}
where $ A(x_0) $ and $ B(x_0) $ are coefficients given by
\begin{eqnarray}
	A(x_0)&=&\Psi\left(\case14+\case{|x_0|}{2}\right)-\Psi\left(\case14\right)+\Psi\left(\case34+\case{|x_0|}{2}\right)-\Psi\left(\case34\right)\label{eq:A}\\
	B(x_0)&=&2^{-1/2-\epsilon}\frac{\Gamma(1/2-\epsilon)}{\Gamma(3/2+\epsilon)}\left[\epsilon\frac{\Gamma(1+\epsilon)}{\Gamma(1-\epsilon)}+\frac{\Gamma(1+\epsilon+|x_0|)}{\Gamma(|x_0|-\epsilon)}\right].\label{B}
\end{eqnarray}
In the first line, $ \Psi(x) $ is the Digamma function, defined by $ \Psi(x)=\case{\rmd}{\rmd x}\log\left(\Gamma(x)\right) $. For $ \epsilon>\frac12 $, i.e. in the positive-recurrent case, the expected first hitting time converges to a constant
\begin{equation}\label{key}
	\langle\tau(0|x_0)\rangle\to\frac{x_0^2}{2\epsilon-1},\qquad\epsilon>\frac12,
\end{equation}
while in the regime $ \epsilon<-\frac12 $, where the walk in transient, we have a probability of hitting the target $ R(0|x_0) <1$, which can be expressed as
\begin{equation}\label{key}
	R(0|x_0)=\frac{(1+\epsilon)_{|x_0|}}{(-\epsilon)_{|x_0|}},\qquad\epsilon<-\frac12,
\end{equation}
where $ (x)_k $ is the Pochhammer symbol \cite{Abr-Steg}:
\begin{equation}\label{key}
	(x)_k=\frac{\Gamma(x+k)}{\Gamma(x)}.
\end{equation}
Let us remark that in the transient regime of the GRW, we have $ R(0|x_0)<1 $ for any choice of $ x_0 $, hence in this case the general correspondence of the first return time probability and the first hitting time probability is always valid also in the transient case, see the discussion at the beginning of section \ref{s:Anal}. 

\subsection{GRW with resetting}
\begin{figure*}[h!]
	\centering
	\begin{tabular}{c @{\quad} c }
		\includegraphics[width=.45\linewidth]{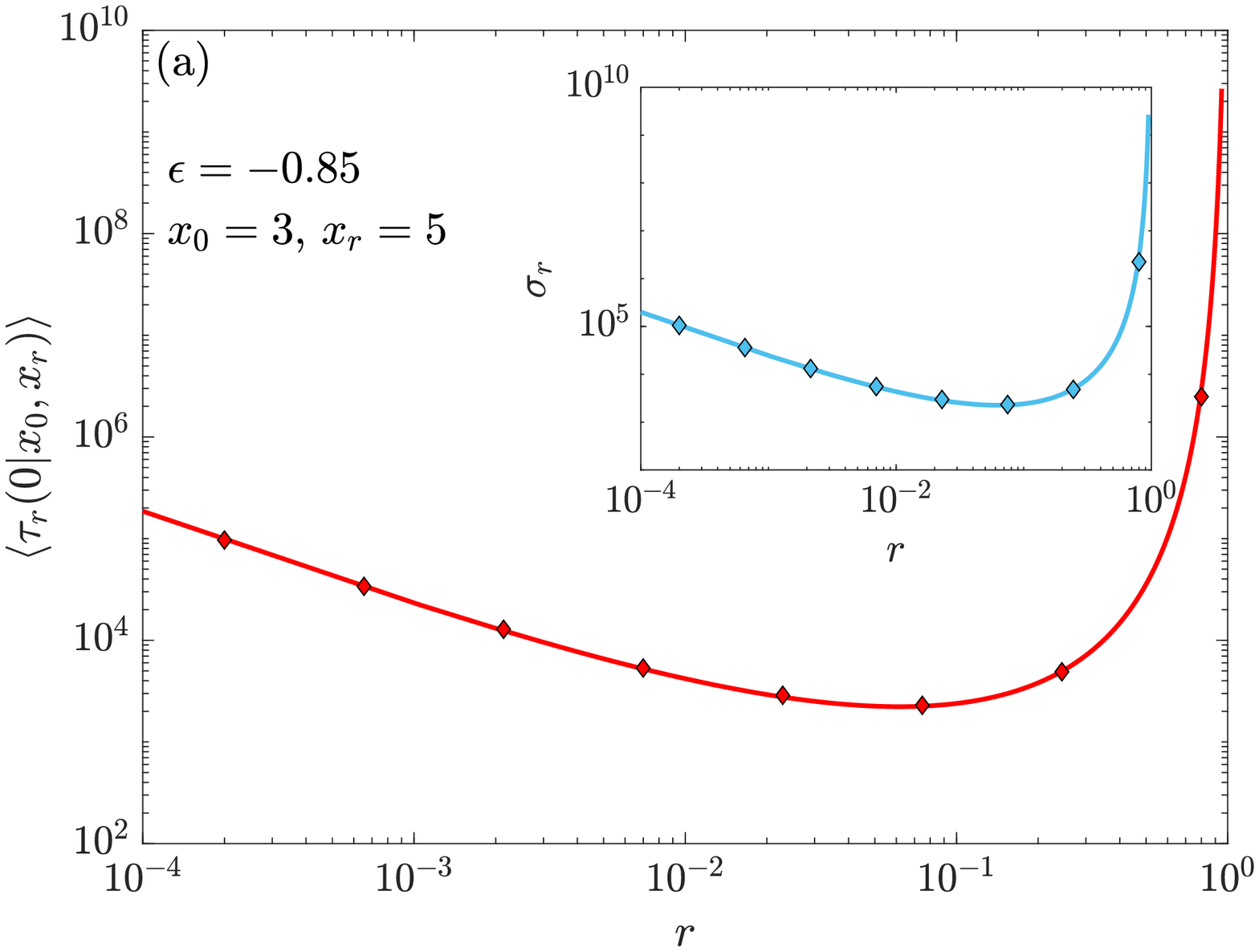} &
		\includegraphics[width=.45\linewidth]{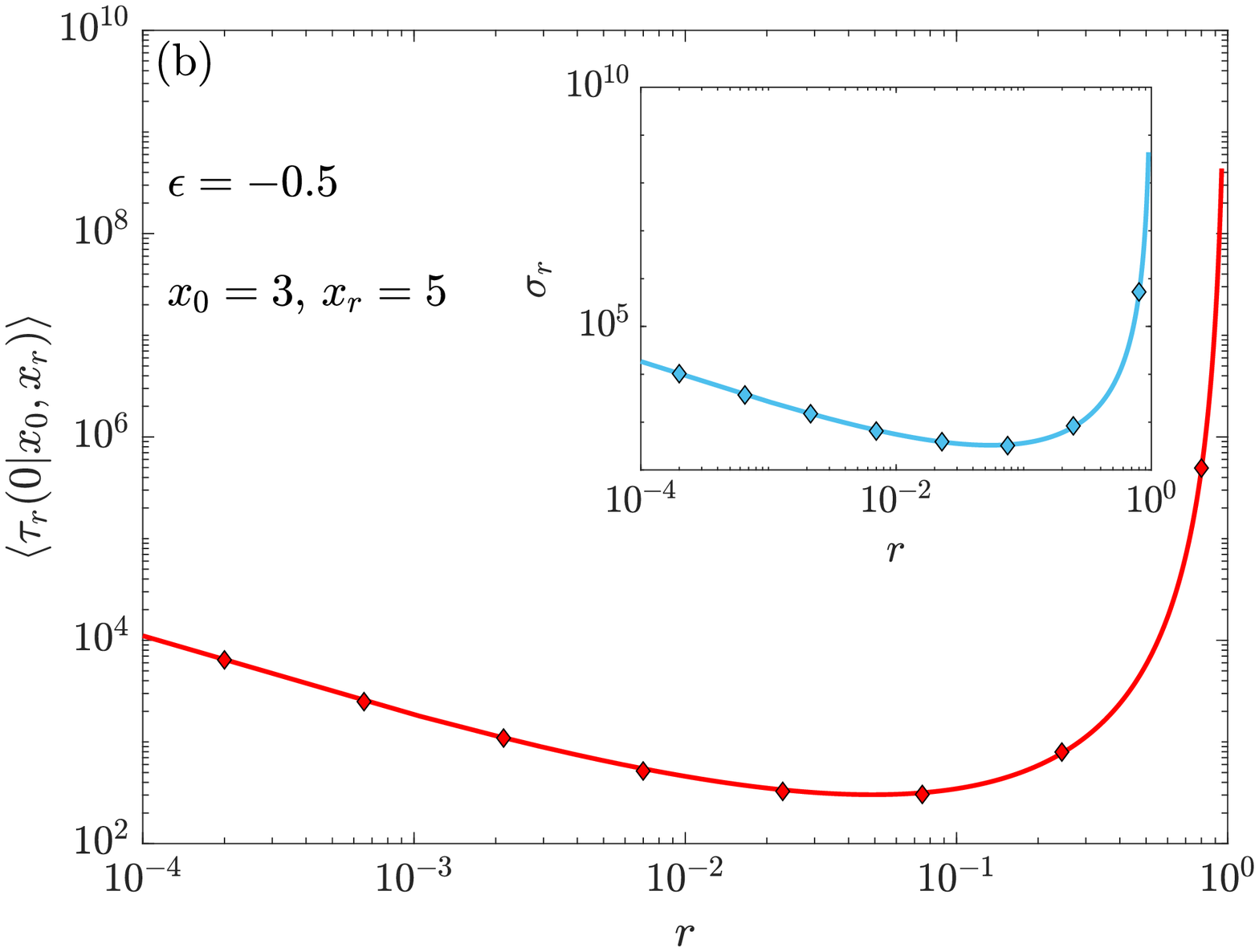} \\
		\includegraphics[width=.45\linewidth]{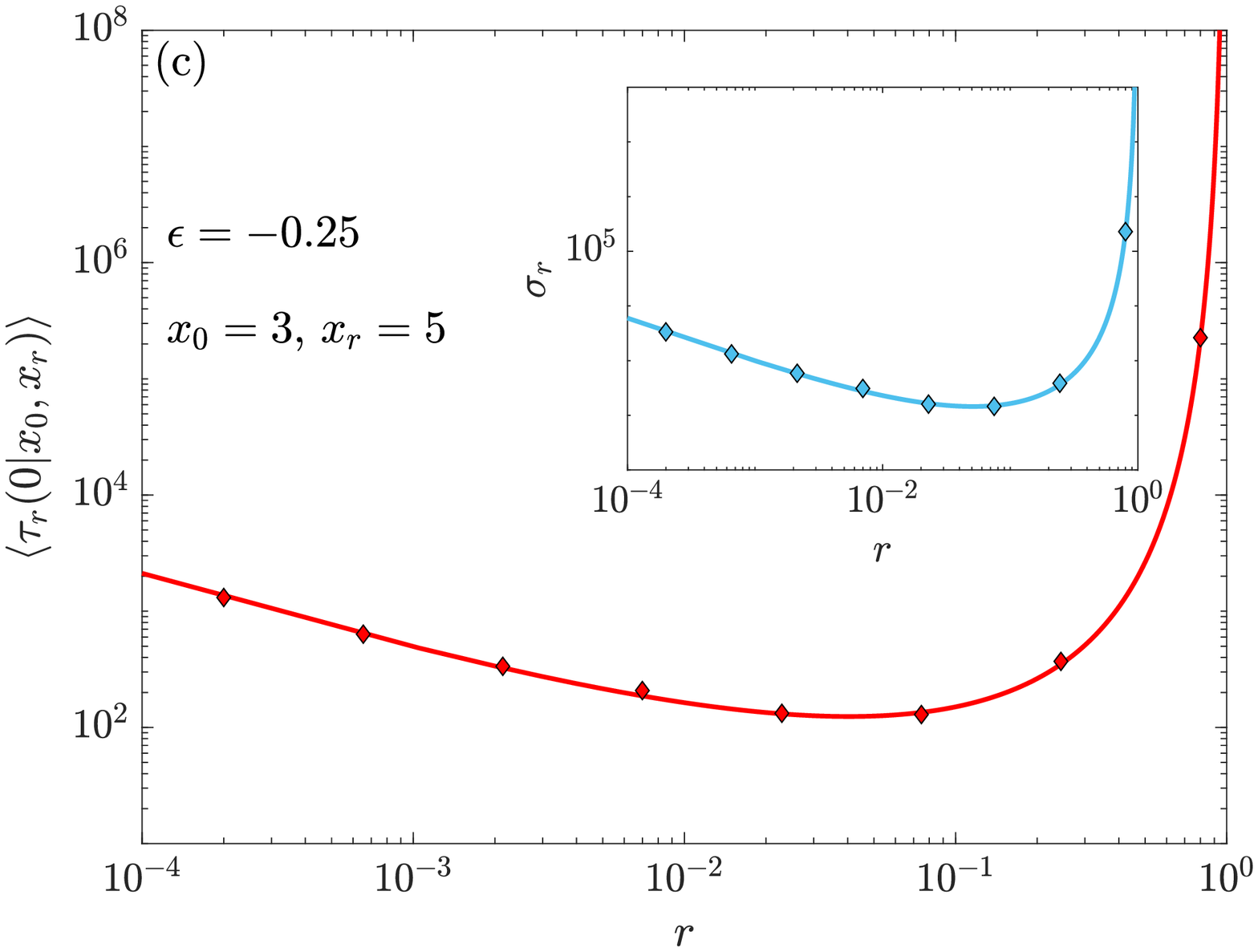} &
		\includegraphics[width=.45\linewidth]{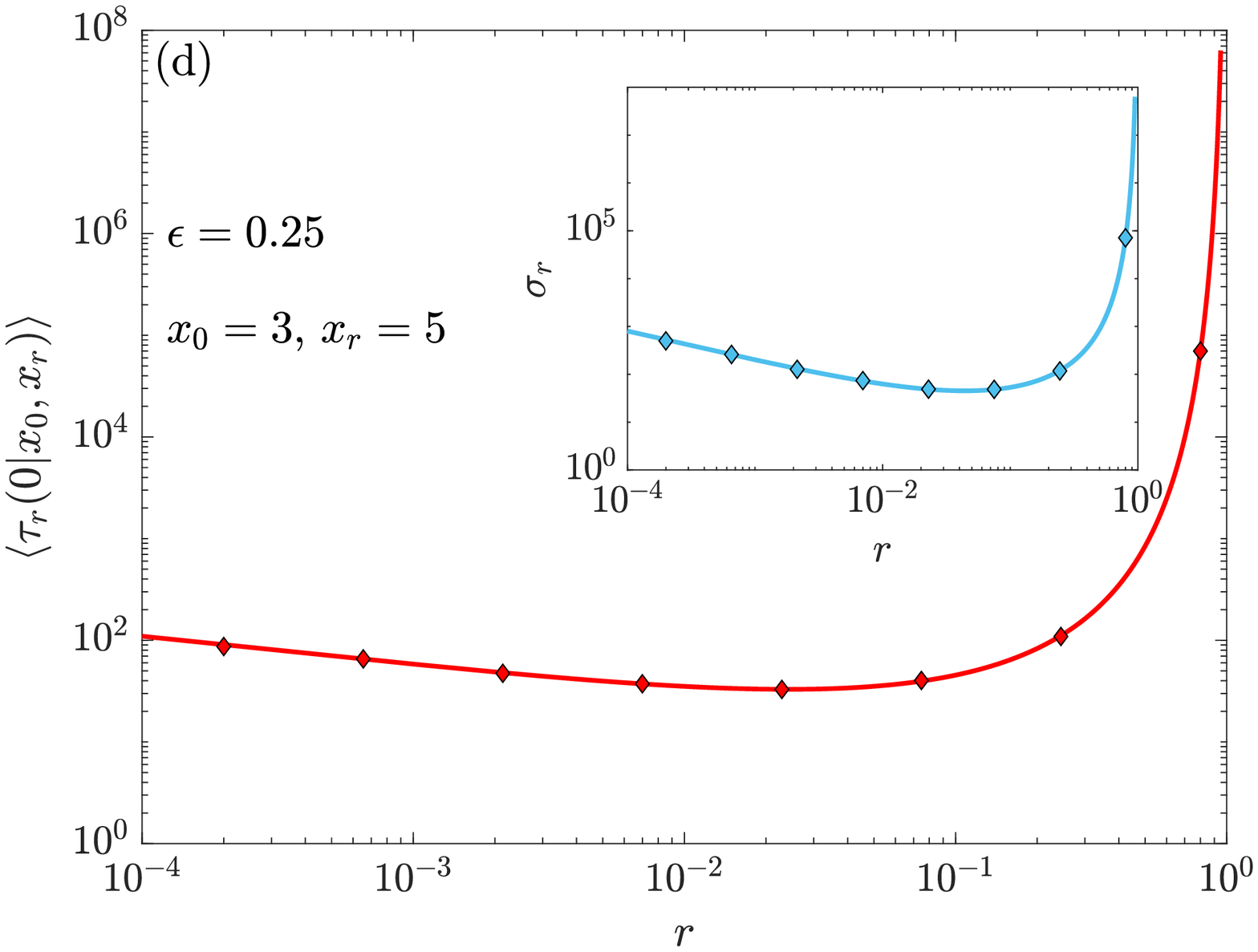}\\
		\includegraphics[width=.45\linewidth]{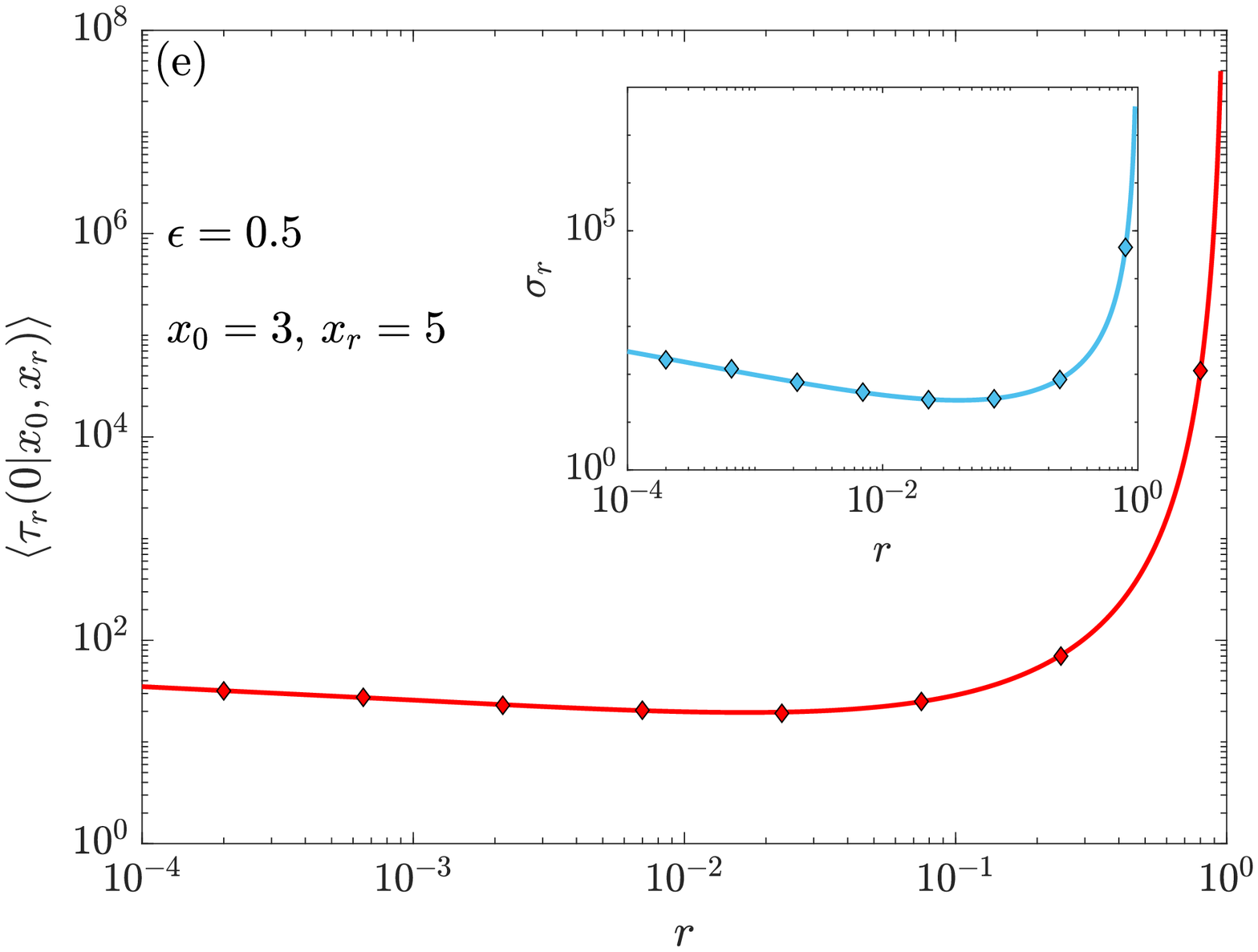} &
		\includegraphics[width=.45\linewidth]{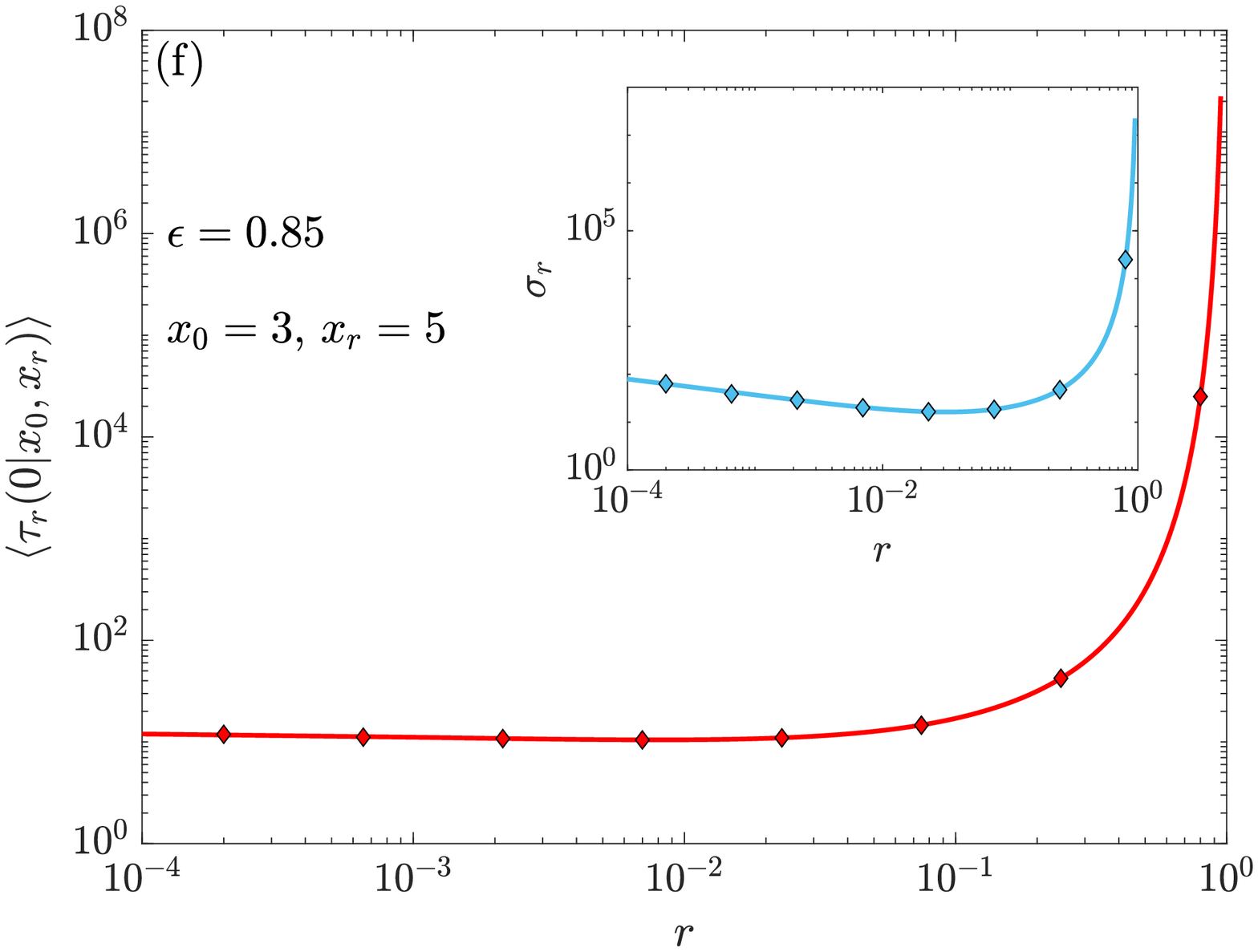}
	\end{tabular}
	\caption{Mean first hitting time and standard deviation (insets) versus the resetting probability $ r $ for different values of $ \epsilon $. Data sets are obtained by averaging over $ N=10^4 $ ($ \epsilon=-0.85,\,-0.5,\,-0.25 $), $ N=10^5 $ ($\epsilon=0.25,\,0.5$) or $ N=10^6 $ ($\epsilon=0.85$) trajectories. The agreement with the theoretical results \eref{eq:tau0r_disp} and \eref{eq:std_dev_disp} is very good.}
	\label{fig:MFHT_sim}
\end{figure*}
We present now the numerical results regarding the first hitting time of the origin for the GRW with resetting. All the simulations have been performed basically by evolving the dynamics according to the transition probability $ G_{xy} $ given by \eref{eq:Trans_def}. We initially fix a value for $ \epsilon $ and set the initial position $ x_0 $ and the resetting position $ x_r $, as well as the resetting probability $ r $. Note that, given the symmetry of the walk, we can restrict ourselves to positive values of $ x_0 $ and $ x_r $. At the first time step, we extract a random number $ u_1 $, uniformly distributed in $ (0,1) $: if $ u_1<r $, then the particle is moved to position $ x_r $ (the resetting can happen even at the first time step); if $ u_1>r $ instead, the dynamics of the GRW takes place, and the particle either moves to $ x_0-1 $, with probability
\begin{equation}\label{key}
	G_{x_0-1,x_0}=\frac12\left(1+\frac{\epsilon}{x_0}\right),
\end{equation}
or to $ x_0+1 $ with the complementary probability. This means that after extracting a second random number $ u_2 $, independent of $ u_1 $ but still uniformly distributed in $ (0,1) $, we check whether $ u_2 $ is smaller or bigger than $G_{x_0-1,x_0}$: in the first case the position is decreased by one unity, otherwise it is increased. The procedure is then iterated until the particle reaches the origin, keeping in mind that the transition probability at each iteration depends on the current position. The total number of iterations corresponds to the first hitting time for that realization, and the first and second moment of the first hitting time distribution are then obtained by averaging over many realizations.

In figure \ref{fig:MFHT_sim} we test the results for the mean and the standard deviation, given by \eref{eq:tau0r_disp} and \eref{eq:std_dev_disp} respectively, for several values of $ \epsilon $, describing the different phases of the GRW. We always consider walks that start from $ x_0=3 $ and are reset to position $ x_r=5 $. Note that in each case the mean first hitting time and the standard deviation diverge when $ r\to0 $ or $ r\to1 $, except when the process is positive-recurrent, which corresponds to panel (f) in the figure. In that case, for $ r\to0 $ the curve of the mean first hitting time attains a finite value, corresponding to mean first hitting time of the underlying process.

We can analytically compute the behaviour of the curves for both $ r\approx1 $ and $ r\approx0 $.  Let us consider first the regime $ r\approx1 $, where we recall that the behaviour of $ \widetilde{F}(0,z|y) $ close to $ z=0 $ is important. From \eref{eq:FGF} we have
\begin{equation}\label{key}
	\widetilde{F}(0,z|y)\sim\frac{z^{|y|}}{|y|!}\frac{\Gamma(1+\epsilon+|y|)}{2^{|y|}\Gamma(1+\epsilon)},
\end{equation}
thus, for any $ -1<\epsilon\leq1 $, we immediately obtain from \eref{eq:tau_asym_r1}
\begin{equation}\label{key}
	\langle\tau_r(0|x_0,x_r)\rangle\sim\frac{|x_r|!\Gamma(1+\epsilon)}{\Gamma(1+\epsilon+|x_r|)}\left(\frac{2}{1-r}\right)^{|x_r|},
\end{equation}
and the same behaviour for the standard deviation, as anticipated by \eref{eq:sig_asym_r1}.

\begin{figure*}[h!]
	\centering
	\begin{tabular}{c @{\quad} c }
		\includegraphics[width=.45\linewidth]{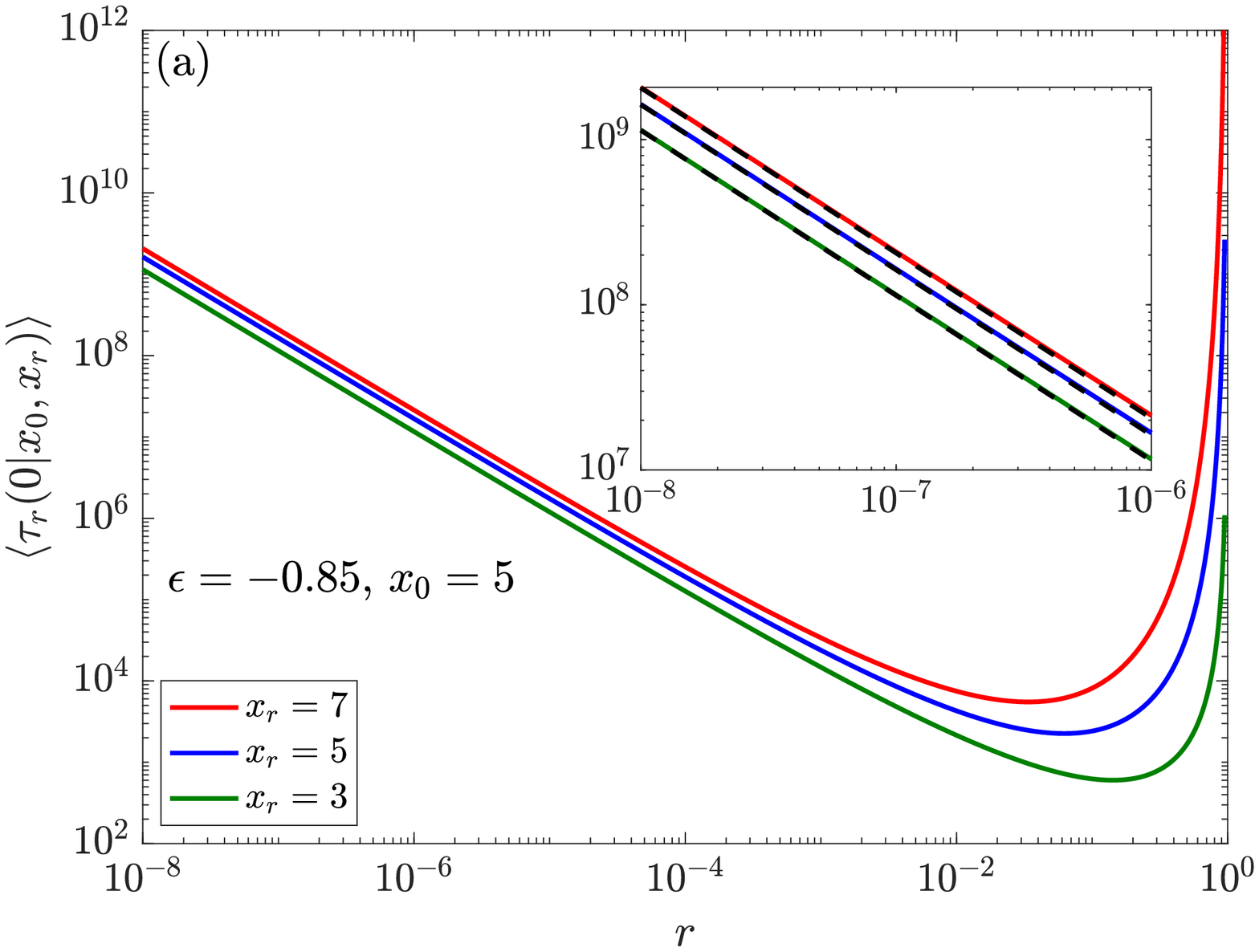} &
		\includegraphics[width=.45\linewidth]{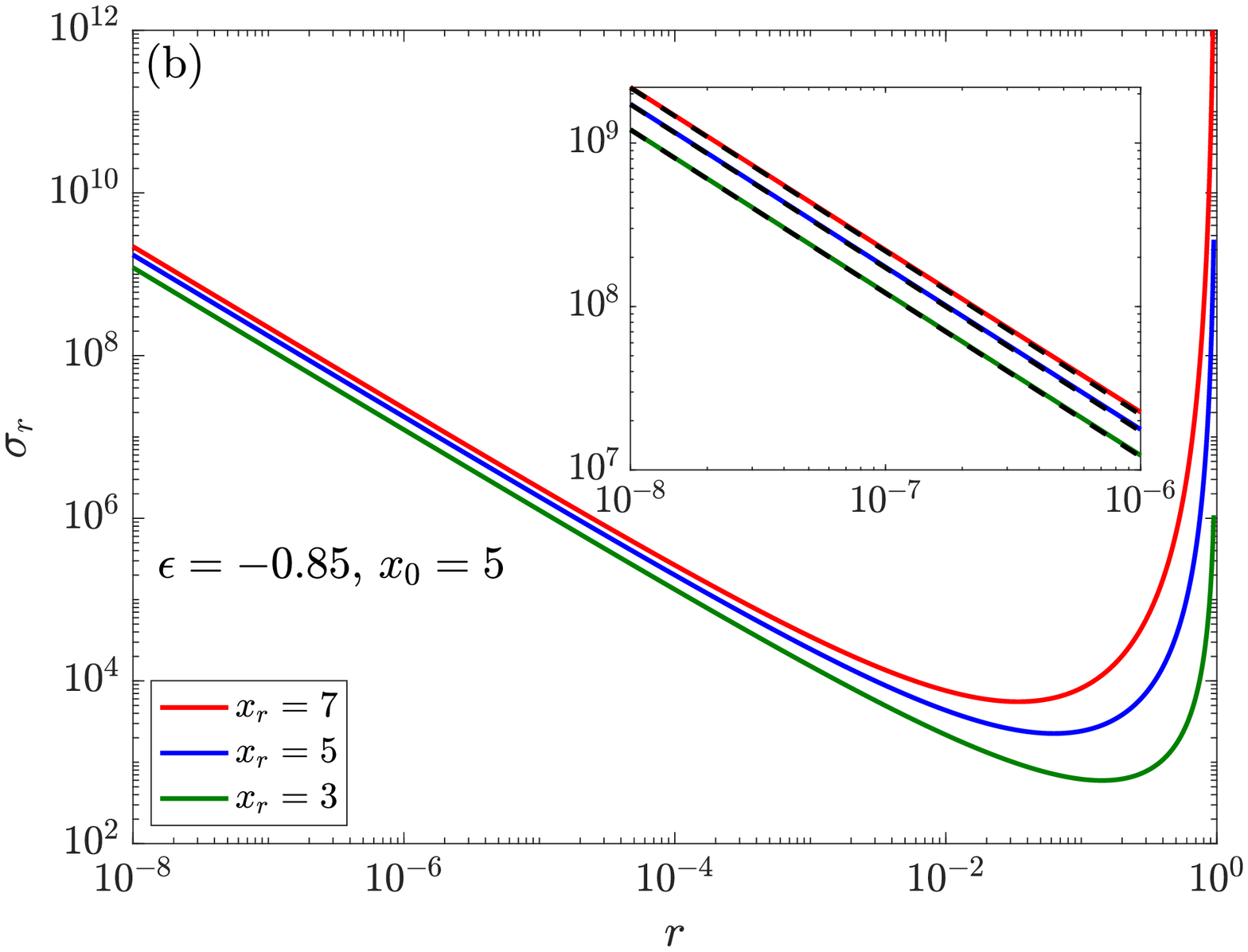}
	\end{tabular}
	\caption{Transient case: mean first hitting time (a) and standard deviation (b) versus the resetting probability $ r $ for different values of the resetting position $ x_r $. The insets show the agreement between the behaviour of the curves for low values of the resetting probability, and the corresponding asymptotics predicted by \eref{eq:MFHT_tr_asy} and \eref{eq:sig_tr_asy}. Note the dependence on $ x_r $.}
	\label{fig:MFHT_tr}
\end{figure*}

Now we take $ r\approx0 $ and begin with the transient case, viz. $ -1<\epsilon<-\case12 $. In this range, for $ z\to1^- $ the generating function $ \widetilde{F}(0,z|y) $ converges to
\begin{equation}\label{key}
	\widetilde{F}(0,z|y)\to R(0|y)=\frac{(1+\epsilon)_{|y|}}{(-\epsilon)_{|y|}},
\end{equation}
hence \eref{eq:tau_asym_r0_tr} yields
\begin{equation}\label{eq:MFHT_tr_asy}
	\langle\tau_r(0|x_0,x_r)\rangle\sim C(x_0,x_r)\frac1r,
\end{equation}
with
\begin{equation}\label{eq:C}
	C(x_0,x_r)=\frac{(-\epsilon)_{x_r}}{(1+\epsilon)_{x_r}}-\frac{\Gamma(1+\epsilon+x_0)\Gamma(x_r-\epsilon)}{\Gamma(1+\epsilon+x_r)\Gamma(x_0-\epsilon)},
\end{equation}
while for the standard deviation we get from \eref{eq:sig_asym_r0_tr}
\begin{equation}\label{eq:sig_tr_asy}
	\langle\sigma_r(0|x_0,x_r)\rangle\sim c(x_0,x_r)\frac1r,
\end{equation}
where
\begin{equation}\label{eq:small_c}
	c(x_0,x_r)=\sqrt{\frac{1-\left[(1+\epsilon)_{x_0}/(-\epsilon)_{x_0}\right]^2}{\left[(1+\epsilon)_{x_r}/(-\epsilon)_{x_r}\right]^2}}.
\end{equation}
A graphic proof is shown in figure \ref{fig:MFHT_tr}, where we observe that the resetting position is relevant to determine the correct asymptotic behaviour for $ r\approx0 $.

\begin{figure*}[h!]
	\centering
	\begin{tabular}{c @{\quad} c }
		\includegraphics[width=.45\linewidth]{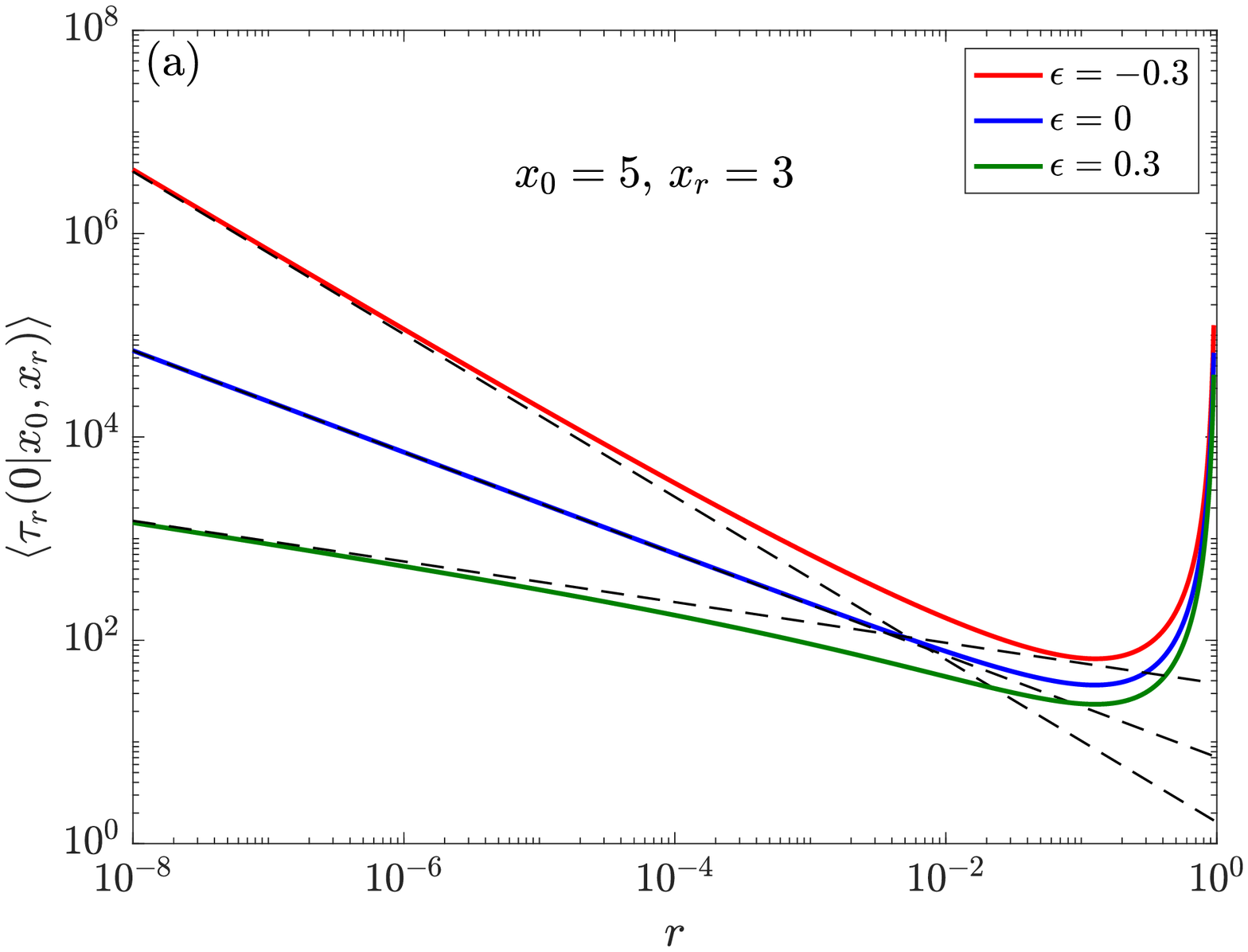} &
		\includegraphics[width=.45\linewidth]{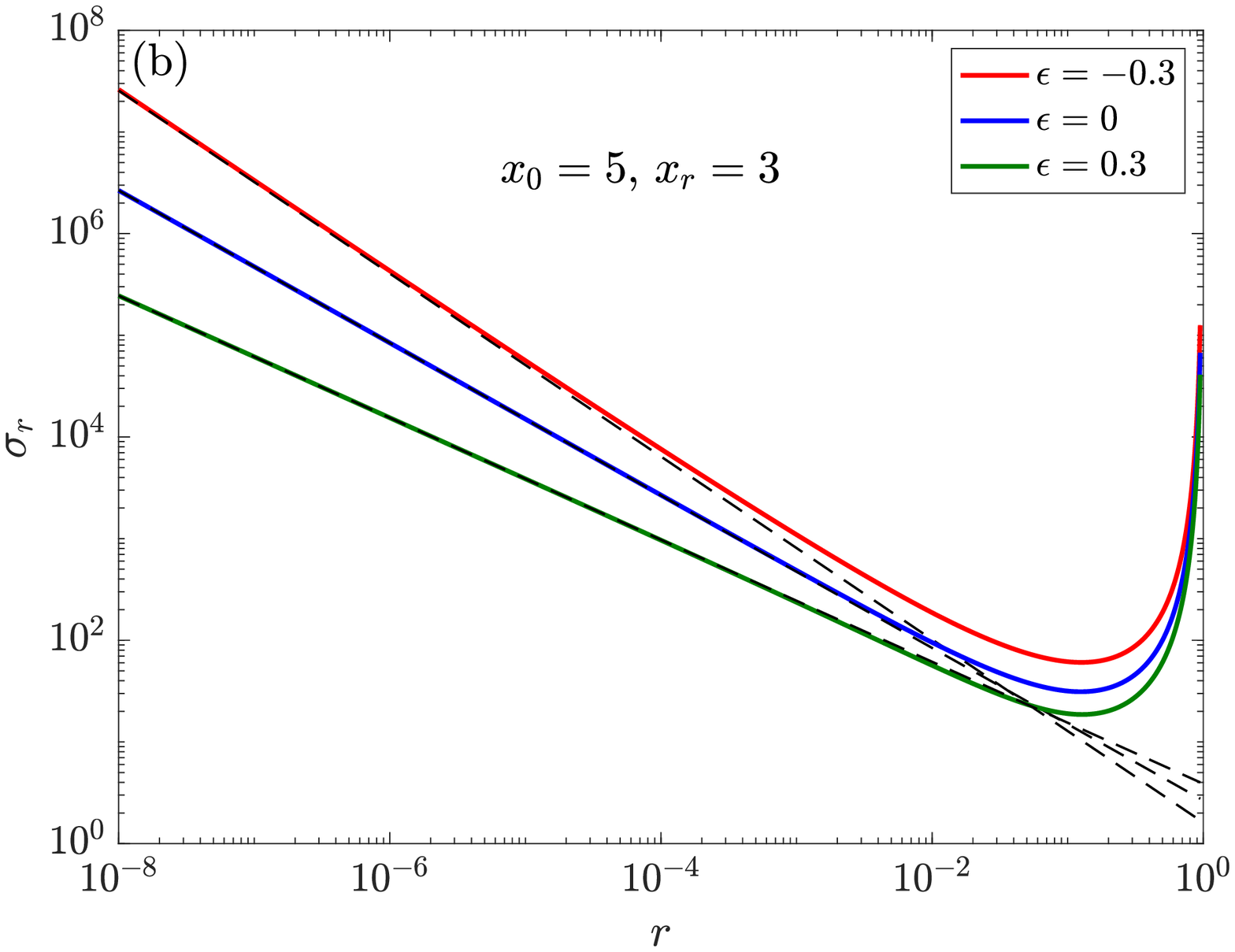}
	\end{tabular}
	\caption{Null-recurrent case: mean first hitting time (a) and standard deviation (b) versus the resetting probability $ r $ for different values of $ \epsilon $. The dashed lines represent the asymptotic behaviour for $ r\approx0 $ predicted by \eref{eq:MFHT_nr_asy} and \eref{eq:sig_nr_asy}.}
	\label{fig:MFHT_nullrec}
\end{figure*}

Now let us consider the transient case and deal with the range $ -\case12\leq\epsilon<\case12 $ first. As anticipated, by using the properties of the hypergeometric function \cite{Abr-Steg} we can write for $ z\approx1 $
\begin{equation}\label{key}
	\widetilde{F}(0,z|y)\sim1-(1-z)^{1/2+\epsilon}L\left(\frac1{1-z};y\right),
\end{equation}
where the slowly-varying function can be determined explicitly, see \ref{app:SV}, obtaining
\begin{equation}\label{key}
	L(t;x)\sim\cases{\frac{A(x)}{\log(t)}&for $\epsilon=-\frac12$\\
	B(x)&for $-\frac12<\epsilon<\frac12$,}
\end{equation}
where $ A(x) $ and $ B(x) $ are the numerical coefficients given by \eref{eq:A} and \eref{B}, respectively. Thus from \eref{eq:tau_asym_r0_rec} we get
\begin{equation}\label{eq:MFHT_nr_asy}
	\langle\tau_r(0|x_0,x_r)\rangle\sim\cases{
		\frac{A(x_0)}{r\log(1/r)}& for $\epsilon=-\frac12$\\
		\frac{B(x_0)}{r^{1/2-\epsilon}}& for $-\frac12<\epsilon<\frac12$,}
\end{equation}
while the standard deviation is obtained from \eref{eq:sig_asym_r0_rec}
\begin{equation}\label{eq:sig_nr_asy}
	\sigma_r\sim\cases{
		\sqrt{\frac{2A(x_0)}{\log(1/r)}}\frac1r& for $\epsilon=-\frac12$\\
		\frac{\sqrt{2\mathcal{B}(x_0)}}{r^{(3/2-\epsilon)/2}}& for $-\frac12<\epsilon<\frac12$,}
\end{equation}
where
\begin{equation}\label{eq:cal_B}
	2\mathcal{B}(x_0)=2^{1/2-\epsilon}\frac{\Gamma(3/2-\epsilon)}{\Gamma(3/2+\epsilon)}\left[\epsilon\frac{\Gamma(1+\epsilon)}{\Gamma(1-\epsilon)}+\frac{\Gamma(1+\epsilon+|x_0|)}{\Gamma(|x_0|-\epsilon)}\right],
\end{equation}
see figure \ref{fig:MFHT_nullrec}.

We now finally deal with the range $ \case12\leq\epsilon\leq1 $, where also the derivative of the slowly-varying function is important. For $ \epsilon=\case12 $, we compute
\begin{equation}\label{key}
	L(t;0,x_0)\sim\frac{x_0^2}{2}\log(t)\Longrightarrow L'(t;0,x_0)\sim\frac{x_0^2}{2t},
\end{equation}
which means we may identify
\begin{equation}\label{key}
	\ell(t;0,x_0)\sim\frac{x_0^2}{2}.
\end{equation}
We therefore obtain that the mean first hitting time grows logarithmically
\begin{equation}\label{key}
	\langle\tau_r(0|x_0,x_r)\rangle\sim\frac{x_0^2}{2}\log\left(1/r\right),
\end{equation}
while the standard deviation diverges as
\begin{equation}\label{eq:sig_rho1_pt}
	\sigma_r\sim\frac{|x_0|}{r^{1/2}}.
\end{equation}
For $ \case12<\epsilon\leq1 $ instead, by expanding $ L(1/(1-z);0,x_0) $ in powers of $ 1-z $ we find
\begin{equation}\label{key}
	L\left(\frac1{1-z};0,x_0\right)\sim \frac{x_0^2}{2\epsilon-1}+B(x_0)(1-z)^{\epsilon-1/2},
\end{equation}
where we point out that $ B(x_0) $ assumes negative values in this case. Thus the mean first hitting time converges to
\begin{equation}\label{key}
	\langle\tau_r(0|x_0,x_r)\rangle\to\frac{x_0^2}{2\epsilon-1},
\end{equation}
and we can say that the derivative of $ L(t;0,x_0) $ is such that
\begin{equation}\label{key}
	tL'(t;0,x_0)\sim-\left(\epsilon-\frac12\right)\frac{B(x_0)}{t^{\epsilon-1/2}}\qquad\textrm{as }t\to\infty.
\end{equation}
By comparing this last equation with \eref{eq:L'_ii}, we find that $ \delta=\epsilon-\case12 $ and
\begin{equation}\label{key}
	\ell(t;0,x_0)=-B(x_0)=\frac1{\delta}\mathcal{B}(x_0),
\end{equation}
where $ \mathcal{B}(x_0) $ is given by \eref{eq:cal_B}, hence the standard deviation diverges as
\begin{equation}\label{eq:sig_rho1_pr}
	\sigma_r\sim\sqrt{2\mathcal{B}(x_0)}r^{-(3/2-\epsilon)/2}.
\end{equation}
However, we point out that in the special case $ \epsilon=1 $ the coefficient $ \mathcal{B}(x_0) $, and hence also $ \sigma_r $, vanish for $ |x_0|=1 $. This is not surprising, because it easily seen from \eref{eq:Trans_def} that with this choice of $ \epsilon $ and $ x_0 $, the walk becomes a trivial deterministic process where the particle always jumps to the origin in one step. The first hitting time distribution has mean one, zero variance, and indeed the generating function given by \eref{eq:FGF_ep1} reduces to
\begin{equation}\label{key}
	\widetilde{F}(0,z|\pm1)=z=1-(1-z).
\end{equation}
An immediate consequence is that resetting can not expedite the process in this case. Apart from physical intuition, this is also confirmed by the criterion discussed in section \ref{s:Opt} related to the coefficient of variation, summarized by \eref{eq:CV_fvar}: if the standard deviation is zero then the lhs can not be bigger than the rhs.

\begin{figure*}[h!]
	\centering
	\begin{tabular}{c @{\quad} c }
		\includegraphics[width=.45\linewidth]{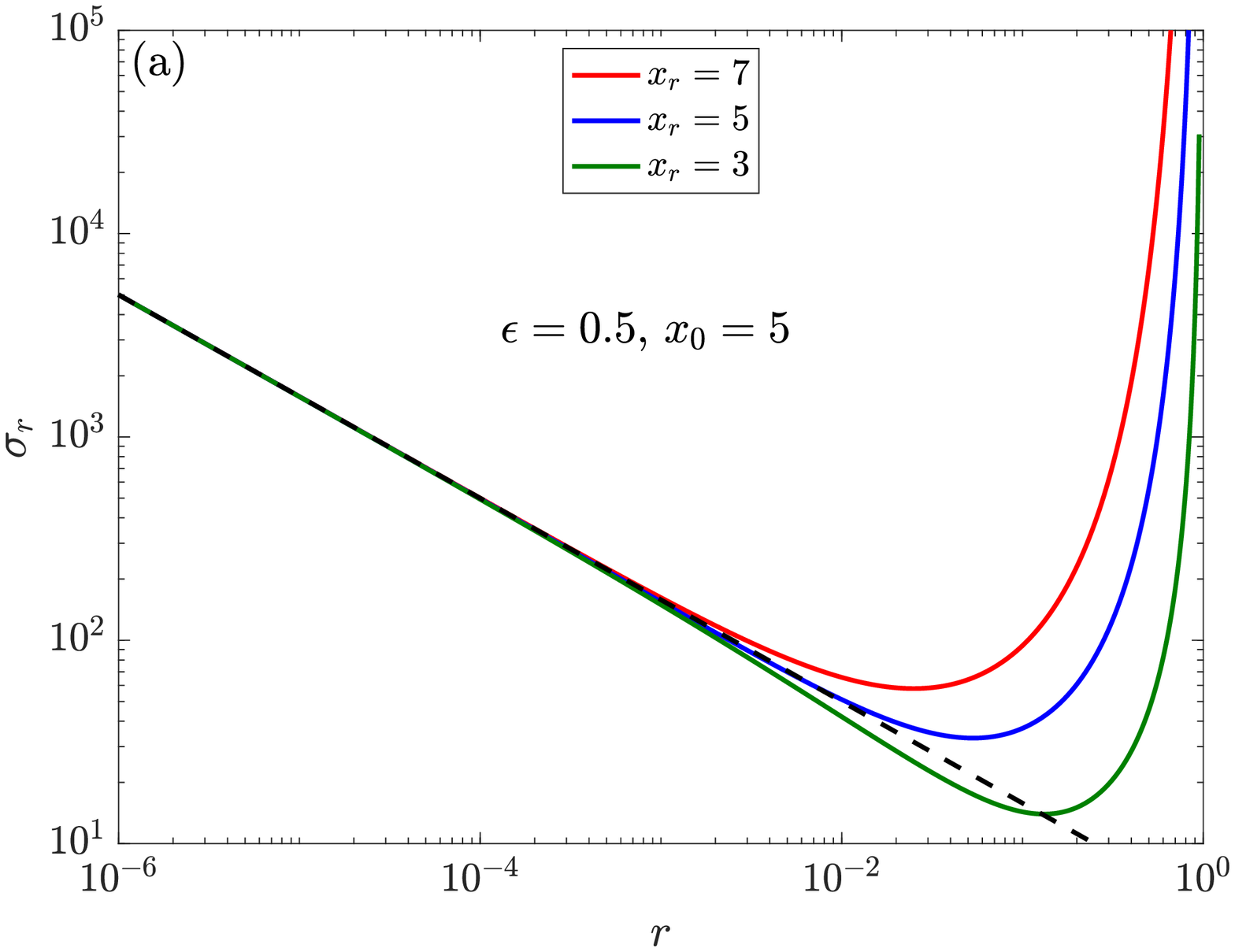} &
		\includegraphics[width=.45\linewidth]{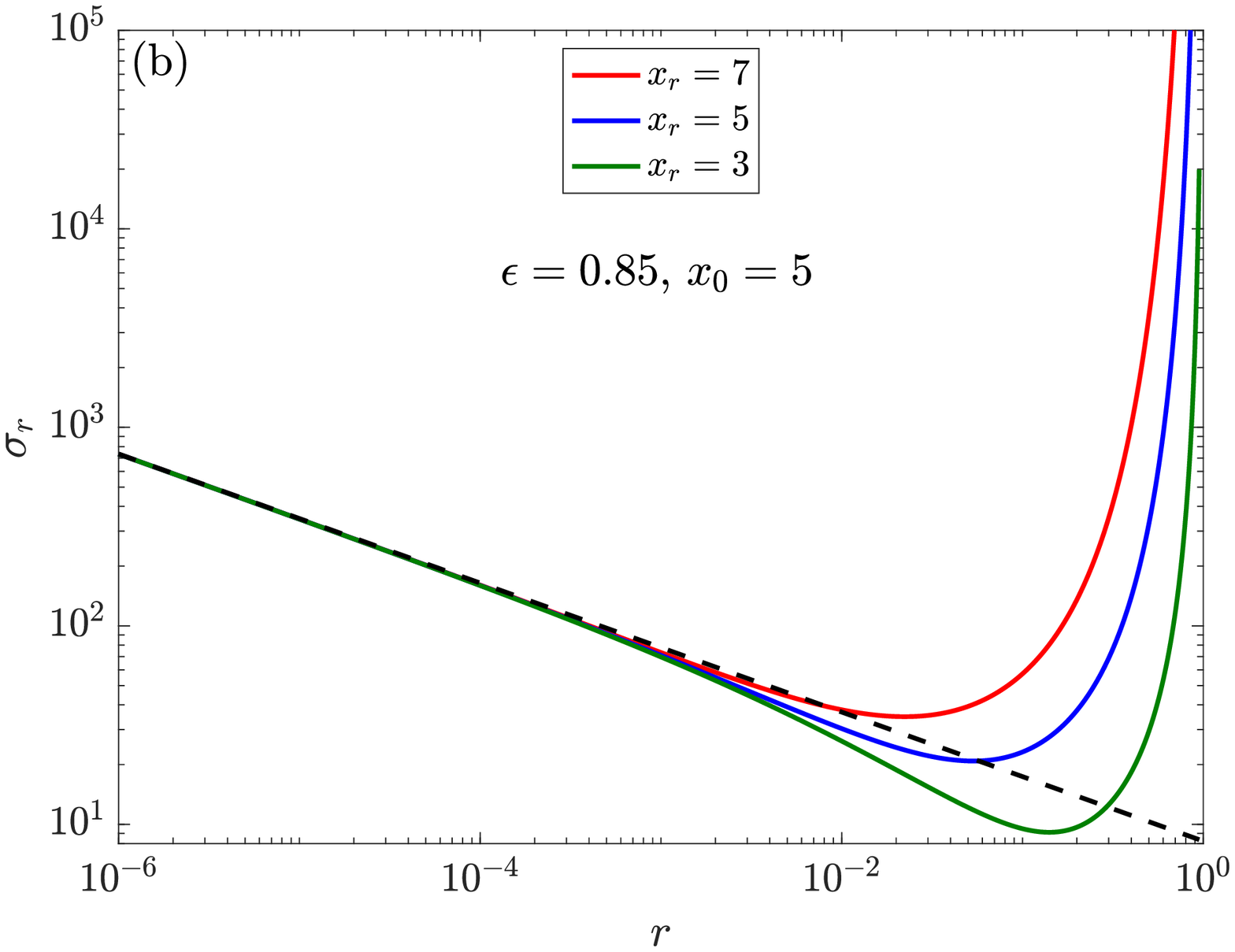}
	\end{tabular}
	\caption{Positive-recurrent case and phase transition: standard deviation versus the resetting probability $ r $ for $ \epsilon=0.5 $ (phase transition) and $ \epsilon=0.85 $ (positive-recurrent case). The dashed lines represent the asymptotic behaviour for $ r\approx0 $ predicted by \eref{eq:sig_rho1_pt} for panel (a) and \eref{eq:sig_rho1_pr} for panel (b).}
	\label{fig:MFHT_rho1}
\end{figure*}

In figure \ref{fig:MFHT_rho1} we present the predicted behaviour of the standard deviation in the range $ \case12\leq\epsilon<1 $, and we observe indeed that as $ r $ becomes small all curves tend to \eref{eq:sig_rho1_pt} for $ \epsilon=\case12 $, and \eref{eq:sig_rho1_pr} for $ \case12<\epsilon<1. $

To summarize, in the small $ r $ regime, the mean first hitting time of the GRW with resetting presents the following behaviour for $ -1<\epsilon\leq\case12 $:
\begin{equation}\label{key}
	\langle\tau_r(0|x_0,x_r)\rangle\sim\cases{C(x_0,x_r)\frac1r & for $-1<\epsilon<-\frac12$\\
		\frac{A(x_0)}{r\log(1/r)}& for $\epsilon=-\frac12$\\
		\frac{B(x_0)}{r^{1/2-\epsilon}}& for $-\frac12<\epsilon<\frac12$\\
		\frac{x_0^2}2\log(1/r)& for $\epsilon=\frac12$,}
\end{equation}
where $ A(x_0) $, $ B(x_0) $ and $ C(x_0,x_r) $ are given by \eref{eq:A}, \eref{B} and \eref{eq:C}, respectively, whereas for $ \epsilon>\case12 $ we recover:
\begin{equation}\label{key}
	\langle\tau_r(0|x_0,x_r)\rangle\to\langle\tau(0|x_0)\rangle=\frac{x_0^2}{2\epsilon-1}.
\end{equation}
In the same regime of $ r $, for the standard deviation we have:
\begin{equation}\label{key}
	\sigma_r\sim\cases{c(x_0,x_r)\frac1r & for $-1<\epsilon<-\frac12$\\
		\sqrt{\frac{2A(x_0)}{\log(1/r)}}\frac1r& for $\epsilon=-\frac12$\\
		\frac{\sqrt{2\mathcal{B}(x_0)}}{r^{(3/2-\epsilon)/2}}& for $-\frac12<\epsilon<\frac12$\\
		\frac{|x_0|}{r^{1/2}}& for $\epsilon=\frac12$\\
		\frac{\sqrt{2\mathcal{B}(x_0)}}{r^{(3/2-\epsilon)/2}}& for $\frac12<\epsilon<1$,}
\end{equation}
with $ c(x_0,x_r) $ given by \eref{eq:small_c}, while $ \mathcal{B}(x_0) $ is defined by \eref{eq:cal_B}. In the limit case $ \epsilon=1 $ we have $ \langle\tau_r(0|x_0,x_r)\rangle\to x_0^2 $, while for $ \sigma_r $:
\begin{equation}\label{key}
	\sigma_r\sim2\sqrt{\frac{|x_0|(x_0^2-1)}{3\sqrt{2r}}}\quad \textrm{for } |x_0|>2,
\end{equation}
where the coefficient results from $ \mathcal{B}(x_0) $ evaluated at $ \epsilon=1 $, and $ \sigma_r\to0 $ for $ |x_0|=1 $.

\subsection{Numerical investigation on the optimal resetting}
We now study numerically the behaviour of the optimal resetting probability for the Gillis model, which we recall is given by the condition
\begin{equation}\label{key}
	\frac{\mathrm{d}}{\mathrm{d}r^*}\langle\tau_{r^*}(0|x_0,x_r)\rangle=0.
\end{equation}
In our case, the solution of this equation shows a dependence on both $ x_0 $ and $ x_r $, as well as $ \epsilon $. However, to reduce the number of free parameters in the problem, we will only deal with the case $ x_r=x_0 $. We recall that in this case we have
\begin{equation}\label{key}
	\langle\tau_{r}(0|x_0,x_0)\rangle=\langle\tau_{r}(0|x_0)\rangle=\frac{1-\widetilde{F}(0,1-r|x_0)}{r\widetilde{F}(0,1-r|x_0)}.
\end{equation}

\begin{figure*}[h!]
	\centering
	\begin{tabular}{c @{\quad} c }
		\includegraphics[width=.45\linewidth]{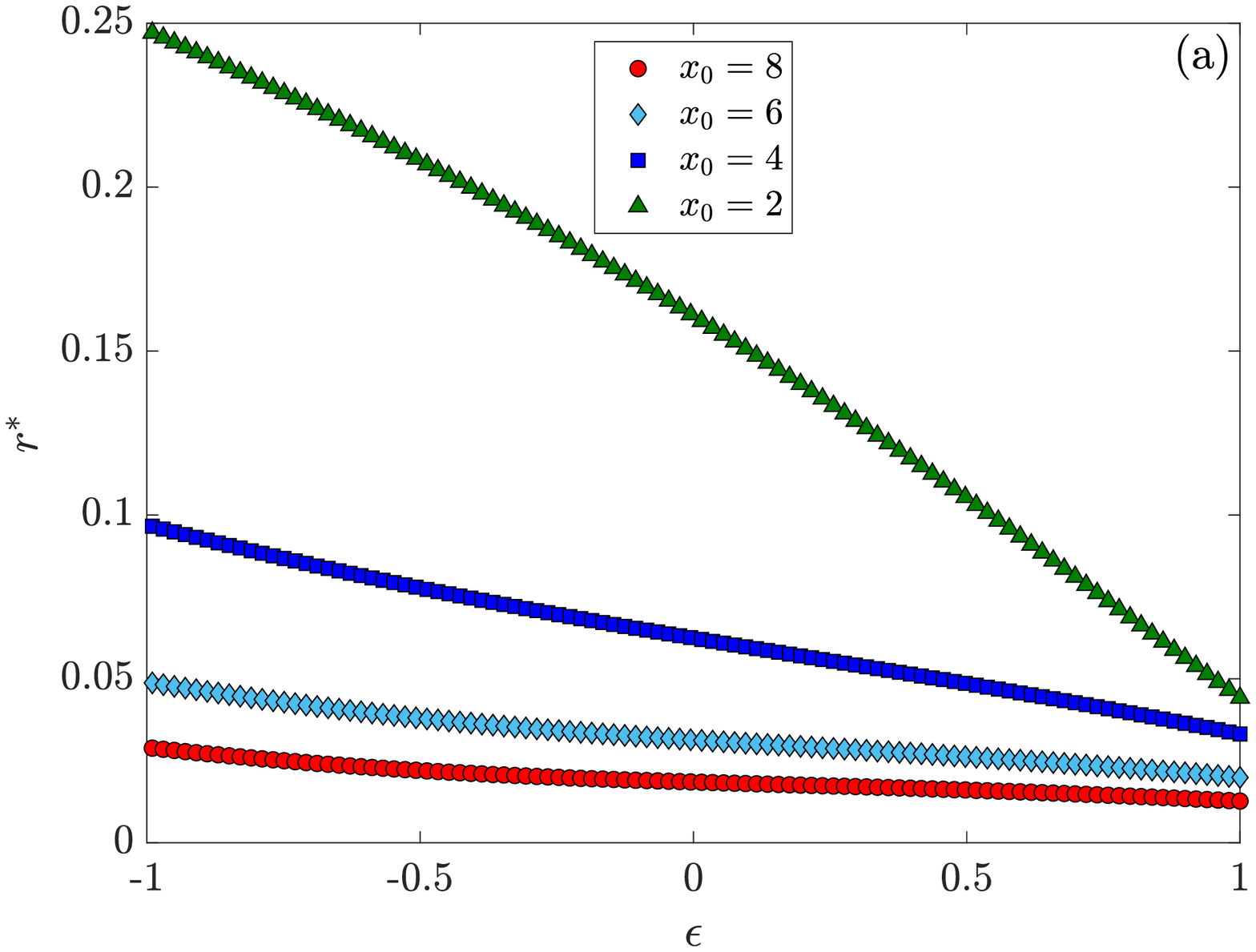} &
		\includegraphics[width=.45\linewidth]{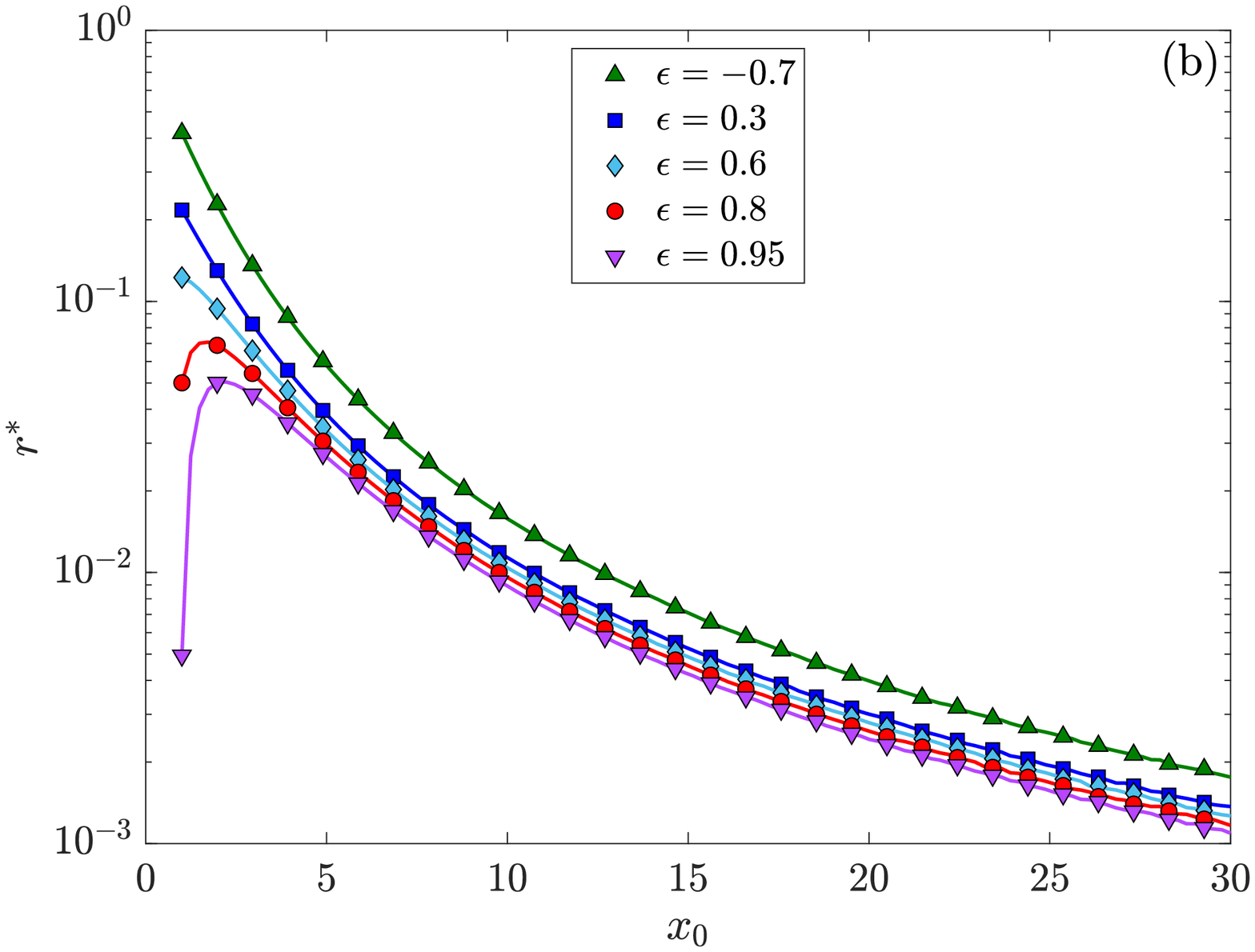}
	\end{tabular}
	\caption{Optimal resetting probability, (a) as a function of $ \epsilon $ for fixed $ x_0 $; (b) as a function of $ x_0 $ for fixed $ \epsilon $. In panel (a), the data converge at both ends of the range to a defined value, which can be obtained for $ \epsilon\approx-1 $ by the solution of \eref{eq:z^*_ep-1}, and for $ \epsilon=1 $ by the solution of \eref{eq:z^*_ep1}}
	\label{fig:ropt_vs}
\end{figure*}

\begin{figure*}[h!]
	\centering
	\begin{tabular}{c @{\quad} c }
		\includegraphics[width=.45\linewidth]{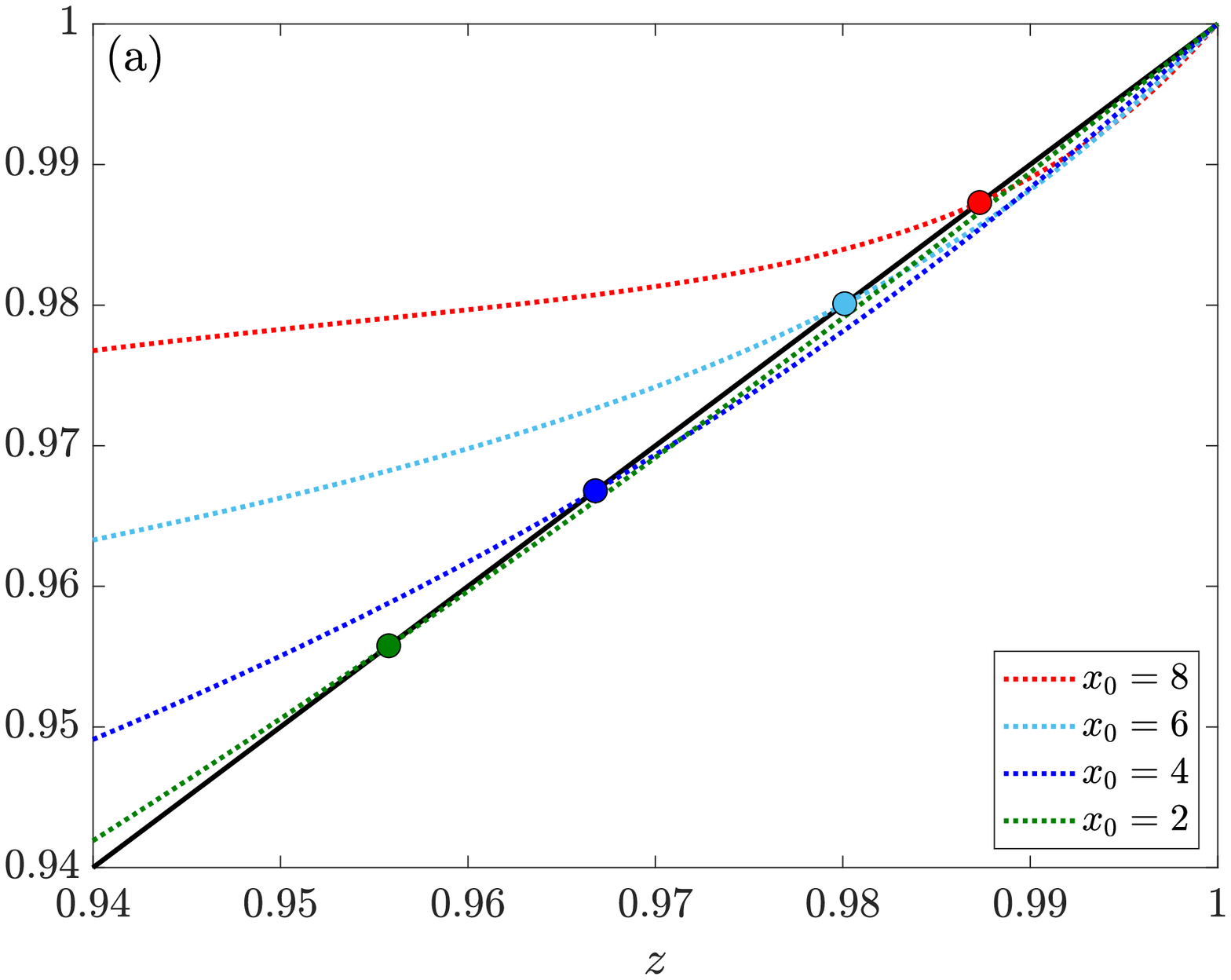} &
		\includegraphics[width=.45\linewidth]{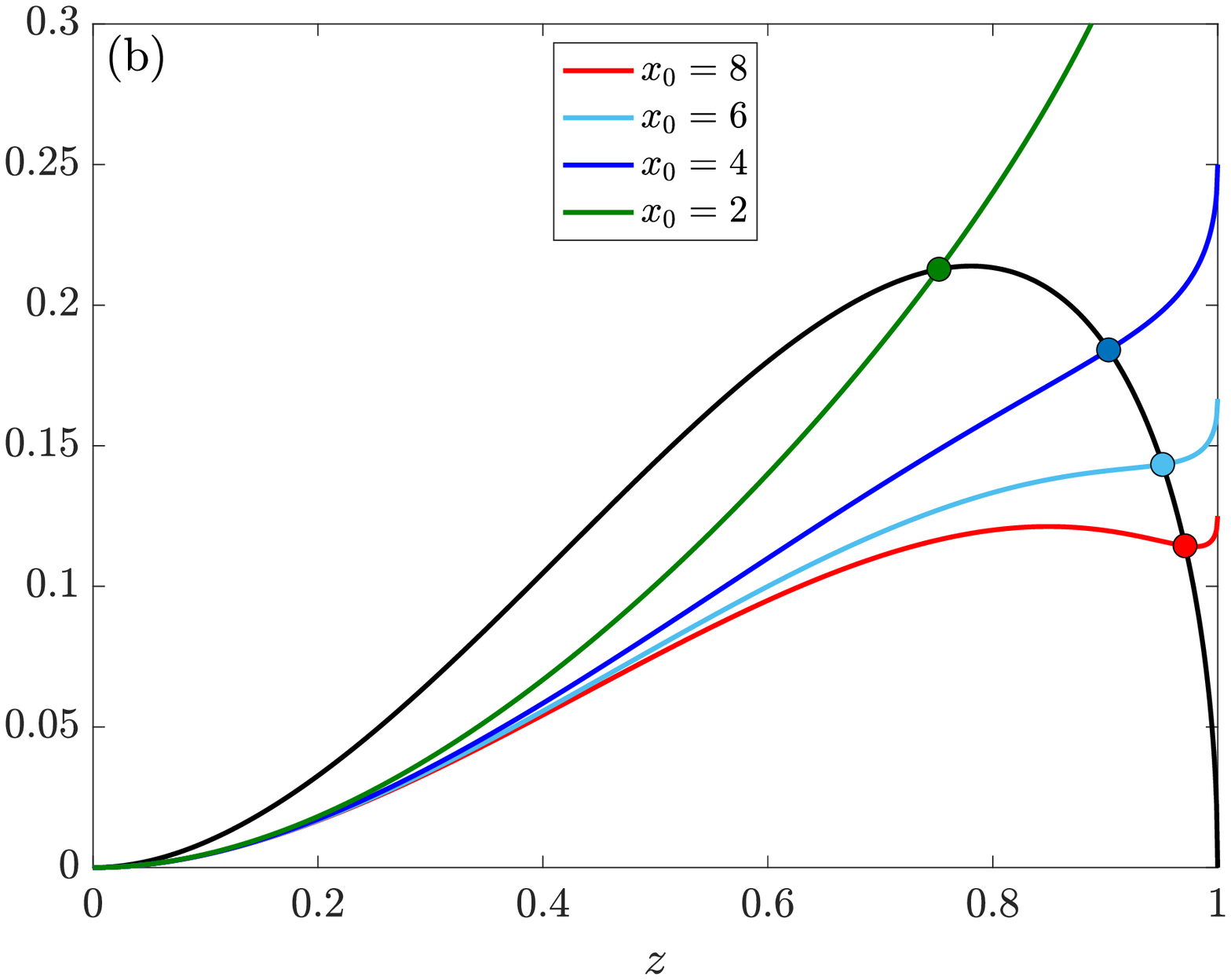}
	\end{tabular}
	\caption{Graphical representation of the solutions of \eref{eq:z^*_ep1} (panel (a)) and  \eref{eq:z^*_ep-1} (panel (b)). The abscissa of each marker corresponds to the optimal value $ z^* $: (a) for $ \epsilon=1 $, we find $ z^*\approx0.9873 $ ($ x_0 =8$), $ 0.9801 $ ($x_0=6$), $ 0.9668 $ ($x_0=4$) and $ 0.9558 $ ($x_0=2$); (b) for $ \epsilon=-1 $, we find $ z^*\approx0.9710 $ ($x_0=8$), $ 0.9610 $ ($x_0=6$), $ 0.9030 $ ($x_0=4$) and $ 0.7522 $ ($x_0=2$).}
	\label{fig:ropt_sol}
\end{figure*}

Figure \ref{fig:ropt_vs} displays the numerically computed values of $ r^* $ and its behaviour as a function of the system parameters. We first note (panel (a)) that the optimal probability shows a monotone decrease with increasing bias $ \epsilon $; furthermore, it converges to a constant value at both ends of the range of existence of $ \epsilon $. To explain this, let us first recall that for $ \epsilon=1 $ the problem is well-defined and can be expressed in terms of the generating function \eref{eq:FGF_ep1}. Then, one can compute $ \langle\tau_{r}(0|x_0)\rangle $ and write an explicit equation for $ r^* $, which can be formulated in terms of $ z^*=1-r^* $ as
\begin{equation}\label{eq:z^*_ep1}
	z^*=\left[h\left(z^*,x_0\right)\right]^{\case1{|x_0|}},
\end{equation}
where
\begin{eqnarray}\label{key}
	h\left(t,x_0\right)&=&\left(1-\sqrt{1-t^2}\right)\left(1+|x_0|\sqrt{1-t^2}\right)-\nonumber\\
	&&|x_0|(1-t)\left[t^{|x_0|-1}-\frac{(|x_0|+1)t^{|x_0|+1}}{(1-\sqrt{1-t^2})(1+|x_0|\sqrt{1-t^2})}\right].
\end{eqnarray}
This can then be solved numerically for $ z^* $, and one finds that the data in figure \ref{fig:ropt_vs} corresponding to $ \epsilon=1 $ coincide indeed with the solution $ r^*=1-z^* $ of the previous equation. The situation close to $ \epsilon=-1 $ is instead less clear, because, as we have already mentioned, the first hitting problem is pathological since the particle can not hit the origin; as a consequence, the generating function \eref{eq:FGF} vanishes for $ \epsilon=-1 $. Nevertheless, observe that for $ \epsilon\approx-1 $  we can write
\begin{equation}\label{key}
	\widetilde{F}(0,z|x_0)\approx K_{\epsilon}(x_0)f(z,x_0),
\end{equation}
where
\begin{equation}\label{key}
	f(z,x_0)=\left(\frac{z}{1+\sqrt{1-z^2}}\right)^{|x_0|}
\end{equation}
is obtained by evaluating the hypergeometric functions appearing in \eref{eq:FGF} at $ \epsilon=-1 $, while
\begin{equation}\label{key}
	K_{\epsilon}(x_0)=\frac{\Gamma(1+\epsilon+|x_0|)}{|x_0|!\,\Gamma(1+\epsilon)},
\end{equation}
is a very small coefficient, due to the divergence of the Gamma function in the denominator. It is easy to verify that the equation for $ r^* $ is
\begin{equation}\label{key}
	\frac{f'(1-r^*,x_0)}{f^2(1-r^*,x_0)}-\frac1{f(1-r^*,x_0)}+K_{\epsilon}(x_0)=0,
\end{equation}
which, since $K_{\epsilon}(x_0)$ is small, can be approximated as
\begin{equation}\label{key}
	\frac{f'(1-r^*,x_0)}{f^2(1-r^*,x_0)}-\frac1{f(1-r^*,x_0)}\approx0.
\end{equation}
This last equation does not depend on $ \epsilon $ and yields a meaningful solution, explaining the behaviour of the data close to $ \epsilon=-1 $ in figure \ref{fig:ropt_vs}. It can be shown that in this case $ r^* $ can be given in terms of $ z^*=1-r^* $, which satisfies
\begin{equation}\label{eq:z^*_ep-1}
	(z^*)^2\sqrt{\frac{1-z^*}{1+z^*}}=\left(1-\sqrt{1-(z^*)^2}\right)\left(1-z^*+\frac{z^*}{|x_0|}\right).
\end{equation}
Let us hence observe that at $\epsilon=-1 $ the system exhibits an abrupt phase transition: although for $ \Delta=\epsilon+1=0 $ there is no $ r $ that can minimise the mean first hitting time, for each positive  and arbitrarily small $ \Delta $, there exists an optimal resetting probability $ r^* $ which is always smaller than $ 1-z^*(-1) $, where $z^*(-1) $ is the solution of \eref{eq:z^*_ep-1}.
From this analysis, and given the monotonic dependence of $ r^* $ on $ \epsilon $ displayed in the figure, we can conclude that for fixed $ x_0 $ the optimal resetting probability will satisfy for any $ -1<\epsilon<1 $
\begin{equation}\label{key}
	1-z^*(1)<r^*<1-z^*(-1),
\end{equation}
where $ z^*(1) $ and $z^*(-1) $ are the solutions of \eref{eq:z^*_ep1} and \eref{eq:z^*_ep-1}, respectively. A graphical illustration of these is offered in figure \ref{fig:ropt_sol}.

It is also interesting to analyse the behaviour of $ r^* $ for fixed $ \epsilon $ as a function of the starting position, panel (b) in figure \ref{fig:ropt_vs}. It is shown that different situations arise depending of the chosen value of the bias. We see that for $ \epsilon \ll 1$, the optimal probability decreases with the distance from the origin. In a random walk setting, where the particle only performs nearest-neighbour jumps, this can be explained by the fact that a longer initial distance requires more steps to reach the target, hence we expect that optimality is achieved by resetting the process less frequently. However, when $ \epsilon $ approaches one, while for $ x_0\geq2 $ the data show the same behaviour, for $ x_0=1 $ the value of $ r^* $ gets closer to zero. This is again due to the fact that for the underlying walk, when $ x_0=1 $ and $ \epsilon=1 $, the variance of the first hitting time vanishes, hence resetting can not expedite the process; thus when $ \epsilon\approx1 $ we also expect the optimal resetting probability to be very small. On the other hand, when $ x_0\geq2 $, the variance of the underlying process is infinite even for $ \epsilon=1 $, hence the data show the monotonic behaviour discussed before.

\begin{figure*}[h!]
	\centering
	\begin{tabular}{c @{\quad} c }
		\includegraphics[width=.45\linewidth]{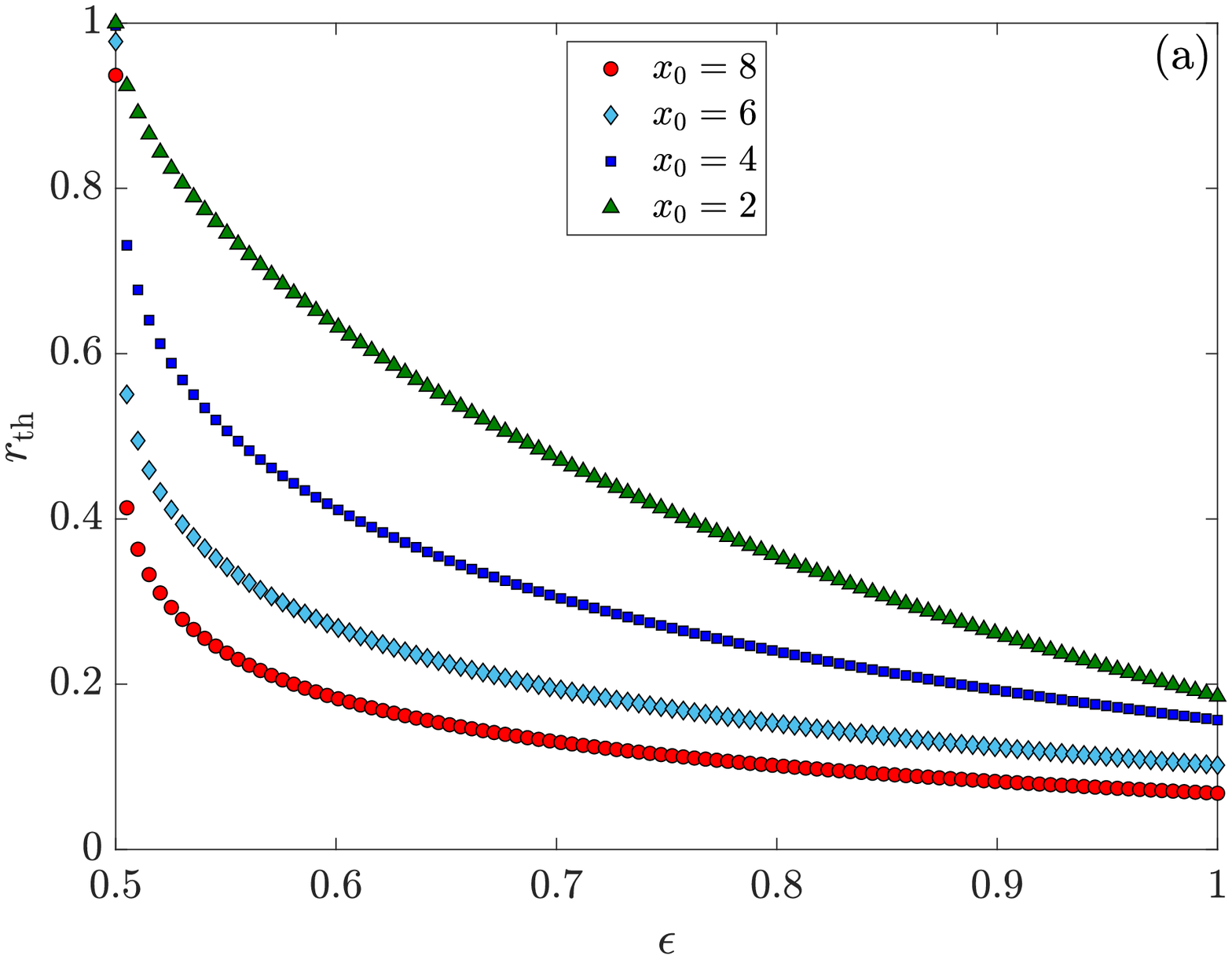} &
		\includegraphics[width=.45\linewidth]{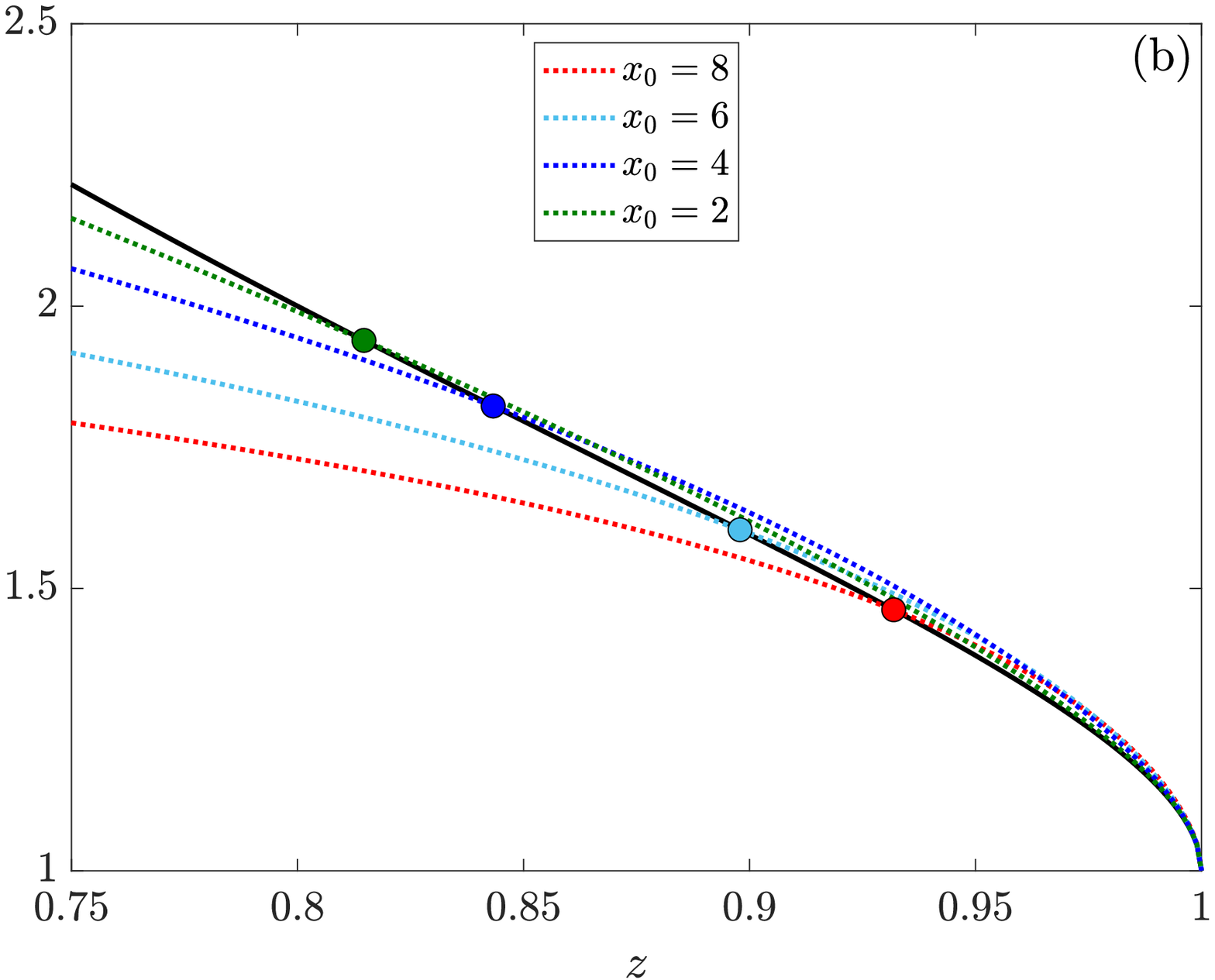}
	\end{tabular}
	\caption{(a) Threshold probability for the positive-recurrent case, as a function of $ \epsilon $. The data approach one for $ \epsilon\approx\case12 $ and reach $ r_{\mathrm{th}}=1-z_{\mathrm{th}} $ for $ \epsilon=1 $, where $ z_{\mathrm{th}} $ is the solution of \eref{eq:rth_ep1}. (b) Graphical representation of the solutions of \eref{eq:rth_ep1}: the abscissa of each marker corresponds to the value of $ z_{\mathrm{th}}=1-r_{\mathrm{th}} $: these are $ z_{\mathrm{th}}\approx0.9319 $ ($ x_0=8 $), $0.8979 $ ($ x_0=6 $), $0.8433 $ ($ x_0=4 $) and $0.8147 $ ($ x_0=2 $).}
	\label{fig:rth}
\end{figure*}

Finally, let us briefly discuss a feature that we can observe only in the positive-recurrent case. As already mentioned, for $ r\to0 $, $ \langle\tau_r(0|x_0)\rangle $ converges to a finite value and $ \rmd \langle\tau_r(0|x_0)\rangle/\rmd r<0 $. On the other hand, for $ r\to1 $ the mean first hitting time diverges, as shown in section \ref{s:Anal}. Therefore, while for sufficiently small $ r $ we have $ \langle\tau_{r}(0|x_0)\rangle <\mu_1(0,x_0)$, there is a threshold value $ r_{\mathrm{th}} $ for the resetting probability above which one observes $ \langle\tau_r(0|x_0)\rangle>\mu_1(0,x_0) $, where $ \mu_1(0,x_0) $ is the mean first hitting time of the reset-free process. This value is hence the solution of
\begin{equation}\label{eq:rth}
	\langle\tau_{r_{\mathrm{th}}}(0|x_0)\rangle=\mu_1(0,x_0).
\end{equation}
A numerical study of the behaviour of $ r_{\mathrm{th}} $ as a function of $ \epsilon $ is presented in figure \ref{fig:rth}, which shows that $ r_{\mathrm{th}} $ does not present particularly surprising features: for $ \epsilon\approx\case12 $, which is the phase transition point, $ r_{\mathrm{th}} $ is close to one, since the underlying process with $ \epsilon=\case12 $ has infinite mean first hitting time; for $ \epsilon\approx1 $ it reaches a finite value, which is the solution of \eref{eq:rth} in the particular case $ \epsilon=1 $. By defining $ z_{\mathrm{th}}=1-r_{\mathrm{th}} $ and using the generating function \eref{eq:FGF_ep1}, it can be shown that this new variable satisfies
\begin{equation}\label{eq:rth_ep1}
	\frac{1+\sqrt{1-z^2_{\mathrm{th}}}}{z_{\mathrm{th}}}=\left\lbrace\left[1+x_0^2\left(1-z^2_{\mathrm{th}}\right)\right]\left(1+|x_0|\sqrt{1-z^2_{\mathrm{th}}}\right)\right\rbrace^{\case1{|x_0|}},
\end{equation}
see figure \ref{fig:rth}. In the intermediate range $ \case12<\epsilon<1 $, the bias towards the origin becomes stronger with increasing $ \epsilon $ and thus the range of $ r $ for which resetting expedites the process becomes smaller, yielding the monotonically decreasing behaviour of $ r_{\mathrm{th}} $ displayed in panel (a).

\section{Conclusions}\label{s:Conc}
We have studied the problem of the first hitting time of the origin for a Gillis random walker starting from $ x_0 $, that can be reset to a generic position $ x_r $ after each step, with fixed probability $ r $. This provides an example of a spatially non-homogeneous random walk for which exact analytical results can be obtained, both in presence and absence of resetting. By starting from a general perspective, we have examined the relation between the recurrence properties of the reset-free system, and the first two moments of the first hitting time distribution of the process with resetting. This analysis let us obtain the exact analytical curves for the mean first hitting time, as well as the variance, which enabled a comparison with numerical simulations. We have also extensively studied the behaviour of the optimal resetting probability for which the mean first hitting time is minimised, as a function of the tuning parameter of the bias $ \epsilon $ and the starting position $ x_0 $, and provided an explanation for the observed behaviour in relation to the physics of the problem.

We stress once again that the system we have considered can be seen as a discrete analogous of a Brownian particle diffusing in an asymptotically logarithmic potential, undergoing stochastic Poissonian resetting. Diffusion with resetting in a logarithmic potential was studied in \cite{RayReu-2020,RayReu-2021}, hence one can notice resemblances between the results presented in this paper and these works. Nevertheless, it has been pointed out in the context of resetting \cite{BonPal-2021,RiaBoyHer-2022} that the discrete-time formalism of random walks presents some important differences with respect to the continuous-time formalism of diffusion processes. For example, in the latter it has been established that the introduction of resetting reduces the mean first passage time whenever the coefficient of variation of the reset-free system is larger than unity \cite{Reu-2016}. For random walks instead this criterion is not valid, as noted in \cite{BonPal-2021,RiaBoyHer-2022} and reported in section \ref{s:Opt}, where we have also extended the result to the case where the starting and the reset position do not coincide. More precisely, the criterion has to be modified to take into account that each step has a fixed duration, and also the resetting requires one unit of time to be performed. Hence, it is no more sufficient that the coefficient of variation is larger than unity, but we can see from \eref{eq:CV_fvar}, by setting $ x_r=x_0 $, that the inequality requires an additional term\textemdash which depends explicitly on the first moment of the underlying first-passage distribution. Similar observations can be made for the condition \eref{eq:CV_opt0} at the optimal reset probability: in the case of diffusion one would simply have that the coefficient of variation evaluated at the optimal reset rate is equal to one, and indeed one can graphically prove this by observing that the abscissa of the intersection point between the curves of mean first passage time and standard deviation corresponds to $ r^* $, see for example \cite{RayReu-2020}. This feature is not observed instead for random walks.

We think that the present work is useful to extend the results on resetting for random walks, by also providing a connection with the problem of Brownian motion in logarithmic potentials, which is of great interest for the physical literature, without neglecting the differences between the discrete-time formalism and the continuous one of diffusion processes.

\appendix
\section{Generating function of the first hitting time for the GRW}\label{a:F_app}
In this Appendix we will compute the generating function $ \widetilde{F}(0,z|x_0) $ in an indirect way, i.e. by making use of the relation \cite{Red}
\begin{equation}\label{eq:F_vs_P}
	\widetilde{F}(0,z|x_0)=\frac{\widetilde{P}(0,z|x_0)}{\widetilde{P}(0,z|0)},
\end{equation}
where $ \widetilde{P}(x,z,|x_0) $ is the generating function of the probability $ P_n(x|x_0) $ of being at site $ x $ at step $ n $ starting from $ x_0 $. A nice method to obtain the function $ \widetilde{P}(0,z|x_0) $ is presented in \cite{Hug-I}, we fully report it here for completeness. We begin by writing the evolution equation for $ P_n(x|x_0) $:
\begin{equation}
	P_{n+1}(x|x_0)=\frac{x-1-\epsilon}{2(x-1)}P_n(x-1|x_0)+\frac{x+1+\epsilon}{2(x+1)}P_n(x+1|x_0).
\end{equation}
By multiplying both sides for $ z^{n+1} $ and summing over $ n=0,1,\dots $, we get an equation involving generating functions which reads
\begin{eqnarray}\label{key}
	\widetilde{P}(x,z|x_0)&=&\delta_{x,x_0}+\frac z2\widetilde{P}(x+1,z|x_0)+\frac z2\widetilde{P}(x-1,z|x_0)\nonumber\\&&+\frac{\epsilon z}{2\left(x+1\right)}\widetilde{P}(x+1,z|x_0)-\frac{\epsilon z}{2\left(x-1\right)}\widetilde{P}(x-1,z|x_0).
\end{eqnarray}
We now consider the discrete Fourier transform of this equation, with
\begin{equation}\label{key}
	\FT(q|x_0)=\sum_{x=-\infty}^{+\infty}\rme^{\rmi qx}\widetilde{P}(x,z|x_0),
\end{equation}
which leads to
\begin{equation}\label{key}
	\frac{1-z\cos q}{\sin q}\FT(q|x_0)-\frac{\rme^{\rmi qx_0}}{\sin q}=-\rmi\epsilon z\sum_{x\neq 0}\frac{e^{iqx}}{x}\widetilde{P}(x,z|x_0).
\end{equation}
By taking the derivative with respect to $ q $ of both sides we obtain a linear first-order differential equation for $ \FT(q|k_0) $, which we put in the following form:
\begin{eqnarray}\label{key}
	\FT'(q|x_0)+f(q)\FT(q|x_0)&=&\frac{\sin q}{1-z\cos q}\frac{\mathrm{d}}{\mathrm{d}q}\left(\frac{\rme^{\rmi qx_0}}{\sin q}\right)\nonumber\\
	&&-\frac{\epsilon z\sin q}{1-z\cos q}\widetilde{P}(0,z|x_0),
\end{eqnarray}
where
\begin{equation}\label{key}
	f(q)=\frac{\sin q}{1-z\cos q}\left[\frac{\mathrm{d}}{\mathrm{d}q}\left(\frac{1-z\cos q}{\sin q}\right)-\epsilon z\right].
\end{equation}
This equation can be solved by means of the usual integrating factor method, and we get the general solution
\begin{eqnarray}\label{key}
	\FT(q|x_0)&=&\frac{\sin q}{\left(1-z\cos q\right)^{1-\epsilon}}\int\frac{\mathrm{d}}{\mathrm{d}q'}\left(\frac{\rme^{\rmi q'x_0}}{\sin q'}\right)\frac{\mathrm{d}q'}{\left(1-z\cos q'\right)^\epsilon}\nonumber\\
	&&-\frac{\epsilon z\sin q\,\widetilde{P}(0,z|x_0)}{\left(1-z\cos q\right)^{1-\epsilon}}\int\frac{\mathrm{d}q'}{\left(1-z\cos q'\right)^\epsilon}+\frac{c\sin q}{\left(1-z\cos q\right)^{1-\epsilon}},
\end{eqnarray}
where we indicated by $ c $ the integration constant. Considering now the integral from $ 0 $ to $ 2\pi $ and dividing both sides by $ 2\pi $ we get, after integrating by parts and rearranging terms:
\begin{equation}\label{key}
	\widetilde{P}(0,z|x_0)=\frac{\int_{0}^{2\pi}\rme^{\rmi qx_0}\left(1-z\cos q\right)^{-1-\epsilon}\mathrm{d}q}{\int_{0}^{2\pi}\left(1-z\cos q\right)^{-\epsilon}\mathrm{d}q}.
\end{equation}
The integrals can be evaluated in terms of hypergeometric functions (see below), yielding
\begin{equation}
	\widetilde{P}(0,z|x_0)=\frac{z^{|x_0|}}{|x_0|!}\frac{\Gamma(1+\epsilon+|x_0|)}{2^{|x_0|}\Gamma(1+\epsilon)}\,\frac{\hyp\left(\frac {1+\epsilon+|x_0|}2,1+\frac {\epsilon+|x_0|}2;1+|x_0|;z^2\right)}{\hyp\left(\frac 12\epsilon,\frac12+\frac 12\epsilon;1;z^2\right)}.
\end{equation}
For $ x_0=0 $ one recovers the classical result by Gillis \cite{Gill-corr}. By plugging this expression in \eref{eq:F_vs_P} we find the result given by \eref{eq:FGF} in the main text. Note that in the case $ x_0=0 $ the relation \eref{eq:F_vs_P} is no longer valid and one must instead use \cite{Red}
\begin{equation}\label{key}
	\widetilde{F}(z)=1-\frac{1}{\widetilde{P}(z)},
\end{equation}
where $ \widetilde{F}(z)=\widetilde{F}(0,z|0) $ and $ \widetilde{P}(z)=\widetilde{P}(0,z|0) $. One then recovers \eref{eq:FGF_Ret} in the main text.

\subsection{Evaluation of the integrals}
Let us consider
\begin{equation}\label{key}
	I=\int_{0}^{2\pi}\frac{\rme^{\rmi qk}}{\left(1-z\cos q\right)^\nu}\mathrm{d}q,
\end{equation}
for $ |z|<1 $, $ k\in\mathbb{Z}$ and $\nu\in\mathbb{R} $. Here $ (\nu)_n $ is the Pochhammer symbol, which can be defined in terms of the Gamma function as
\begin{equation}\label{key}
	(\nu)_n=\frac{\Gamma(\nu+n)}{\Gamma(\nu)}.
\end{equation}
An expansion of the integrand in powers of $ z $ yields
\begin{equation}\label{key}
	I=\sum_{n=0}^{\infty}\frac{(\nu)_n}{n!}z^n\int_{0}^{2\pi}\rme^{\rmi qk}\cos^nq\mathrm{d}q.
\end{equation}
For each value of $ n $, by using the binomial expansion of $ \cos^n q $, one deduces that the integral in the $ n $-th term of the series vanishes unless $ n\geq |k| $, and $ n $ and $ k $ are both even or odd, in which case one gets:
\begin{equation}\label{key}
	\int_{0}^{2\pi}\rme^{\rmi kq}\cos^nq\mathrm{d}q=\frac{2\pi}{2^n}{n\choose \frac{n-k}{2}}.
\end{equation}
Let us first consider the case $ n $ even, so that it can be replaced by $ 2m $. The integral is then transformed into the series
\begin{equation}\label{app:hyp:eq:Series_1}
	I=2\pi\sum_{m=\frac{|k|}{2}}^{\infty}\frac{(\nu)_{2m}}{2^{2m}}\frac{z^{2m}}{\Gamma\left(1+m-\frac{k}{2}\right)\Gamma\left(1+m+\frac{k}{2}\right)},
\end{equation}
or, by the simple shift of the index $ m\to m+\frac {|k|}2 $:
\begin{equation}\label{app:hyp:eq:Series_2}
	I=\frac{2\pi z^{|k|}}{2^{|k|}}\sum_{m=0}^{\infty}\frac{(\nu)_{2m+|k|}}{2^{2m}}\frac{z^{2m}}{\Gamma(m+1)\Gamma(1+m+|k|)}.
\end{equation}
For $ n $ odd, we set $ n=2m+1 $ and obtain a series similar to \eref{app:hyp:eq:Series_1}, but starting from $ m=\frac {|k|-1}2 $ and with $ 2m $ replaced by $ 2m+1 $ in each term. Nevertheless, by shifting the index $ m\to m+\frac{|k|-1}2 $, we obtain the same series of \eref{app:hyp:eq:Series_2}. From the definition of the Pochhammer symbol, we deduce that
\begin{equation}\label{key}
	(\nu)_{2m+|k|}=(\nu)_{|k|}(\nu+|k|)_{2m}=\frac{\Gamma(\nu+|k|)}{\Gamma(\nu)}(\nu+|k|)_{2m},
\end{equation}
and moreover from the properties of the Gamma function, we have
\begin{equation}\label{key}
	(\nu+|k|)_{2m}=2^{2m}\left(\frac{\nu+|k|} 2\right)_m\left(\frac{1+\nu+|k|}{2}\right)_m.
\end{equation}
Plugging the last two formulas in Eq. \eref{app:hyp:eq:Series_2}, we finally obtain:
\begin{equation}\label{key}
	I=\frac{2\pi}{|k|!}\left(\frac z2\right)^{|k|}\frac{\Gamma(\nu+|k|)}{\Gamma(\nu)}\,\hyp\left(\frac{\nu+|k|}{2},\frac{1+\nu+|k|}{2};1+|k|;z^2\right).
\end{equation}

\section{Generating function of the first hitting time for the GRW: evaluation of the derivative when $ \rho=1 $ }\label{a:Fr_der}
In this Appendix we are interested in the behaviour of the derivative $ \widetilde{F}'(x,z|y) $ when $ z\approx1 $, knowing that
\begin{equation}\label{eq:app_Fstr}
	\widetilde{F}(x,z|y)\sim1-(1-z)L\left(\frac1{1-z};x,y\right),
\end{equation}
where $ L(t;x,y) $ is a slowly-varying function of $ t $. To simplify the notation, from now on we will not indicate the dependence on the parameters $ x $ and $ y $. We begin by writing
\begin{equation}\label{eq:app_F'}
	\widetilde{F}'(x,z|y)\sim L\left(\frac1{1-z}\right)-\frac{1}{1-z}L'\left(\frac1{1-z}\right),
\end{equation}
where, as we have already pointed out in section \ref{s:Anal2}, the dominant term is the first one at the rhs. However, to correctly capture the subleading correction to the asymptotic behaviour, now we are interested also in the second one, which can be evaluated by making some considerations on $ \widetilde{F}(z) $. Depending on the statistical properties of the underlying first hitting process, we can distinguish three cases:
\begin{enumerate}
\item If the first hitting time distribution of the underlying process has finite mean $ \mu_1 $ and finite second moment $ \mu_2 $, while \eref{eq:app_Fstr} holds, it must also hold
\begin{equation}\label{key}
	\widetilde{F}(z)\sim1-(1-z)\mu_1+\frac{(1-z)^2}{2}(\mu_2-\mu_1).
\end{equation}
This can be seen by expanding $ \widetilde{F}(z) $ in powers of $ 1-z $. It follows that
\begin{equation}\label{key}
	L\left(\frac1{1-z}\right)\sim\mu_1-\frac{1-z}{2}(\mu_2-\mu_1),
\end{equation}
viz.
\begin{equation}\label{key}
	L(t)\sim\mu_1-\frac{1}{2t}(\mu_2-\mu_1)
\end{equation}
for $ t\to\infty $, and we can thus say that
\begin{eqnarray}
	L(t)&\to&\mu_1\\
	L'(t)&\sim&\frac{1}{2t^2}\left(\mu_2-\mu_1\right).
\end{eqnarray}

\item If the first hitting time distribution of the underlying process has finite mean $ \mu_1 $, but infinite second moment, i.e. the first hitting time probability decays as $ F(n)\propto n^{-2-\delta} $, with $ 0<\delta<1 $, then we can write
\begin{equation}\label{key}
	\widetilde{F}(z)\sim1-(1-z)\mu_1+(1-z)^{1+\delta}\ell\left(\frac1{1-z}\right),
\end{equation}
where $ \ell(t) $ is a slowly-varying function. We may thus identify
\begin{equation}\label{key}
	L\left(\frac1{1-z}\right)\sim\mu_1-(1-z)^{\delta}\ell\left(\frac1{1-z}\right),
\end{equation}
hence for $ t\to\infty $:
\begin{eqnarray}
	L(t)&\to&\mu_1\\
	L'(t)&\sim&\delta\frac{\ell(t)}{t^{1+\delta}}.
\end{eqnarray}

\item If we examine the limit case $ \delta=0 $, then also the mean diverges. In particular, it diverges as the slowly-varying function $ L(t) $ (typically displaying a logarithmic growth, because here $ F(n)\propto n^{-2 }$). Then in this case we may write
\begin{equation}\label{key}
	\widetilde{F}'(z)\sim L\left(\frac1{1-z}\right)-\ell\left(\frac1{1-z}\right),
\end{equation}
where $ L(t) $ dominates $ \ell(t) $ as $ t\to\infty $, and so we can identify
\begin{equation}\label{key}
	L'\left(\frac1{1-z}\right)\sim(1-z)\ell\left(\frac1{1-z}\right),
\end{equation}
showing that $ L'(t) $ decays as $ \ell(t)/t $. In general, since we typically expect to observe a logarithmic growth of the mean first hitting time, we also expect that $ L'(t) $ converges to zero as $ 1/t $; it follows from the previous equation that $ \ell(t) $ is expected to converge to a constant. 

\end{enumerate}

\section{Evaluation of the slowly-varying functions}\label{app:SV}
In this Appendix we evaluate, for $ -\case12\leq\epsilon\leq1 $, the slowly-varying functions that appear in the expression of $ \widetilde{F}(0,z|x_0) $ for the GRW:
\begin{equation}\label{eq:app_FGF}
	\widetilde{F}(0,z|x_0)=\frac{z^{|x_0|}}{|x_0|!}\frac{\Gamma(1+\epsilon+|x_0|)}{2^{|x_0|}\Gamma(1+\epsilon)}\,\frac{\hyp\left(\frac {1+\epsilon+|x_0|}2,1+\frac {\epsilon+|x_0|}2;1+|x_0|;z^2\right)}{\hyp\left(\frac {1+\epsilon}2,1+\frac 12\epsilon;1;z^2\right)}.
\end{equation}

For $ \epsilon\neq\pm\case12 $, we use the following property of the hypergeometric functions \cite{Abr-Steg}:
\begin{eqnarray}\label{eq:app_pr1}
	\hyp\left(a,b;c;z\right)&=&(1-z)^{c-a-b}\frac{\Gamma(c)\Gamma(a+b-c)}{\Gamma(a)\Gamma(b)}\Bigg[G(1-z)+\nonumber\\
	&&(1-z)^{a+b-c}\frac{\Gamma(c-a-b)\Gamma(a)\Gamma(b)H(1-z)}{\Gamma(a+b-c)\Gamma(c-a)\Gamma(c-b)}\Bigg],
\end{eqnarray}
where
\begin{eqnarray}
	G(z)&=&\hyp\left(c-a,c-b;c-a-b+1;z\right)\\
	H(z)&=&\hyp\left(a,b;a+b-c+1;z\right).	
\end{eqnarray}
By plugging \eref{eq:app_pr1} in \eref{eq:app_FGF} we can write
\begin{equation}\label{key}
	\widetilde{F}(0,z|x_0)=z^{|x_0|}\frac{G_N(1-z^2)+q_N(1-z^2)^{1/2+\epsilon}H_N(1-z^2)}{G_D(1-z^2)+q_D(1-z^2)^{1/2+\epsilon}H_D(1-z^2)},	
\end{equation}
with
\begin{eqnarray}
	q_N&=&-2^{-1-2\epsilon}\frac{\Gamma(1/2-\epsilon)\Gamma(1+\epsilon+|x_0|)}{\Gamma(3/2+\epsilon)\Gamma(|x_0|-\epsilon)}\\
	q_D&=&\frac{\epsilon}{2^{1+2\epsilon}}\frac{\Gamma(1/2-\epsilon)\Gamma(1+\epsilon)}{\Gamma(3/2+\epsilon)\Gamma(1-\epsilon)}.	
\end{eqnarray}
We now consider an expansion in powers of $ 1-z $ for $ z\to1^- $. When $ -\case12<\epsilon<0 $ we obtain
\begin{eqnarray}\label{key}
	\widetilde{F}(0,z|x_0)&\sim&1-(1-z)^{1/2+\epsilon}\bigg[2^{1/2+\epsilon}\left(q_D-q_N\right)\nonumber\\
	&&-2^{1+2\epsilon}q_D\left(q_D-q_N\right)(1-z)^{1/2+\epsilon}+O\left((1-z)^{1/2-\epsilon}\right)\bigg],	
\end{eqnarray}
hence by writing the coefficients explicitly, we have
\begin{eqnarray}\label{key}
	L\left(\frac1{1-z};x_0\right)&=& B(x_0)-\frac{\epsilon B(x_0)}{2^{1/2+\epsilon}}\frac{\Gamma(1/2-\epsilon)\Gamma(1+\epsilon)}{\Gamma(3/2+\epsilon)\Gamma(1-\epsilon)}(1-z)^{1/2+\epsilon}\nonumber\\
	&&+o\left((1-z)^{1/2+\epsilon}\right),
\end{eqnarray}
where $ B(x_0) $ is given by \eref{B} in the main text; when $ 0<\epsilon<\case12 $ we obtain instead
\begin{eqnarray}\label{key}
	\widetilde{F}(0,z|x_0)&&\sim1-(1-z)^{1/2+\epsilon}\bigg[2^{1/2+\epsilon}\left(q_D-q_N\right)\nonumber\\
	&&-\left(2g_N-2g_D-|x_0|\right)(1-z)^{1/2-\epsilon}+O\left((1-z)^{1/2+\epsilon}\right)\bigg],	
\end{eqnarray}
with
\begin{eqnarray}
	g_N&=&\frac12\frac{(1+|x_0|-\epsilon)(|x_0|-\epsilon)}{1-2\epsilon}\\
	g_D&=&\frac\epsilon2\frac{1-\epsilon}{2\epsilon-1}.	
\end{eqnarray}
It follows
\begin{equation}\label{key}
	L\left(\frac1{1-z};x_0\right)= B(x_0)-\frac{x_0^2}{1-2\epsilon}(1-z)^{1/2-\epsilon}+o\left((1-z)^{1/2-\epsilon}\right).
\end{equation}
Finally when $ \case12<\epsilon\leq1 $ the expansion is
\begin{eqnarray}\label{key}
	\widetilde{F}(0,z|x_0)&&\sim1-(1-z)\bigg[|x_0|+2g_D-2g_N\nonumber\\
	&&+2^{1/2+\epsilon}\left(q_D-q_N\right)(1-z)^{\epsilon-1/2}+O\left((1-z)^{1/2+\epsilon}\right)\bigg],	
\end{eqnarray}
hence
\begin{equation}\label{key}
	L\left(\frac1{1-z};x_0\right)= \frac{x_0^2}{2\epsilon-1}+B(x_0)(1-z)^{\epsilon-1/2}+o\left((1-z)^{\epsilon-1/2}\right).
\end{equation}

We now consider the critical cases $ \epsilon=\pm\case12 $, when \eref{eq:app_pr1} is not valid. For $ \epsilon=-\case12 $ we shall use \cite{Abr-Steg}:
\begin{eqnarray}\label{key}
	\hyp\left(a,b,;a+b;z\right)&=&\frac{\Gamma(a+b)}{\Gamma(a)\Gamma(b)}\log\left(\frac1{1-z}\right)\times\nonumber\\
	&&\left[\sum_{n=0}^{\infty}\frac{(a)_n(b)_n}{(n!)^2}(1-z)^n\left(1-\frac{p(n)}{\log\left(1-z\right)}\right)\right],
\end{eqnarray}
with $ p(n) $ defined by
\begin{equation}\label{key}
	p(n)=2\Psi(n+1)-\Psi(a+n)-\Psi(b+n).
\end{equation}
We can thus write
\begin{equation}\label{key}
	\widetilde{F}(0,z|x_0)=z^{|x_0|}\frac{\sum_{n=0}^{\infty}k_N(n)(1-z^2)^n\left[1-\frac{p_N(n)}{\log(1-z^2)}\right]}{\sum_{n=0}^{\infty}k_D(n)(1-z^2)^n\left[1-\frac{p_D(n)}{\log(1-z^2)}\right]},
\end{equation}
and by considering an expansion in powers of $ 1-z $, we obtain
\begin{eqnarray}\label{key}
	\widetilde{F}(0,z|x_0)&&\sim1-\frac{1}{\log(1-z)}\bigg[p_N(0)-p_D(0)+\nonumber\\
	&&(2k_N(1)-2k_D(1)-|x_0|)(1-z)\log\left(\frac1{1-z}\right)+O(1-z)\bigg]
\end{eqnarray}
hence
\begin{equation}\label{key}
	L\left(\frac1{1-z};x_0\right)=\frac{A(x_0)}{\log\left(\frac1{1-z}\right)}-\frac{x_0^2}{2}(1-z)+o(1-z),
\end{equation}
where $ A(x_0) $ is given by \eref{eq:A} in the main text.

For $ \epsilon=\case12 $ we need the formula \cite{Abr-Steg}
\begin{eqnarray}\label{key}
	\fl\hyp\left(a,b;a+b-1;z\right)=&&\frac{\Gamma(a+b-1)}{\Gamma(a)\Gamma(b)}(1-z)^{-1}\Bigg\{1+(a-1)(b-1)\times\nonumber\\
	&&\sum_{n=0}^{\infty}\frac{(a)_n(b)_n}{n!(n+1)!}(1-z)^{n+1}\bigg[\log(1-z)+\Psi(a+n)+\nonumber\\
	&&\Psi(b+n)-\Psi(n+1)-\Psi(n+2)\bigg]\Bigg\},
\end{eqnarray}
whereby we write
\begin{equation}\label{key}
	\widetilde{F}(0,z|x_0)=z^{|x_0|}\frac{1+c_N\sum_{n=0}^{\infty}r_N(n)(1-z^2)^{n+1}\left[\log(1-z)+s_N(n)\right]}{1+c_D\sum_{n=0}^{\infty}r_D(n)(1-z^2)^{n+1}\left[\log(1-z)+s_D(n)\right]}
\end{equation}
from which it follows
\begin{eqnarray}\label{key}
	\widetilde{F}(0,z|x_0)&\sim&1-(1-z)\bigg[2(c_N-c_D)\log\left(\frac1{1-z}\right)-2c_Ns_N(0)+\nonumber\\
	&&2c_Ds_D(0)+|x_0|+O((1-z)\log(1-z))\bigg].
\end{eqnarray}
By writing the coefficients explicitly, we have
\begin{equation}\label{key}
	L\left(\frac1{1-z};x_0\right)=\frac{x_0^2}{2}\log\left(\frac1{1-z}\right)-K(x_0)+o(1),
\end{equation}
where $ K(x_0) $ is a constant that can be determined from the previous formulae.

\section*{References}
\providecommand{\newblock}{}

\end{document}